\documentclass[12pt]{article}
\usepackage[german, french, american]{babel}
\usepackage[latin1]{inputenc}
\usepackage[T1]{fontenc}
\usepackage{amsmath, amssymb, amsfonts}
\usepackage{amsthm}
\usepackage{mathtools}
\usepackage{mathrsfs}
\usepackage{enumerate}
\usepackage{graphicx}
\usepackage{color}
\usepackage{marvosym}
\usepackage{booktabs}

\usepackage{fancyhdr}

\newtheorem{theorem}{Theorem}[section]
\newtheorem{lemma}[theorem]{Lemma}
\newtheorem{remark}[theorem]{Remark}
\newtheorem{corollary}[theorem]{Corollary}

\numberwithin{equation}{section}

\newcommand{\MA}{\text{\Male}}
\newcommand{\VE}{\text{\Female}}

\newcommand{\RR}{\mathbb{R}}
\newcommand{\NN}{\mathbb{N}}

\numberwithin{equation}{section}

\usepackage{fix-cm}

\begin{document}

\title{Analysis and Numerics for an Age- and Sex-Structured Population Model}

\author{
Michael Pokojovy\thanks{Department of Mathematics and Statistics, University of Konstanz, Konstanz, Germany \hfill 
\texttt{michael.pokojovy@uni-konstanz.de}} \and
Yevhenii Skvarkovskyi\thanks{Department of Cybernetics, Kyiv National Taras Shevchenko University, Ukraine \hfill 
\texttt{y.skvarkovsky@gmail.com}}}

\date{\today}

\maketitle

\pagestyle{fancy}
\thispagestyle{plain}

\begin{abstract}
	We study a linear model of McKendrick-von Foerster-Keyfitz type for the temporal development of the age structure of a two-sex human population.
	For the underlying system of partial integro-differential equations,
	we exploit the semigroup theory to show the classical well-posedness and asymptotic stability in a Hilbert space framework
	under appropriate conditions on the age-specific mortality and fertility moduli.
	Finally, we propose an implicit finite difference scheme to numerically solve this problem
	and prove its convergence under minimal regularity assumptions.
	A real data application is also given.
\end{abstract}


{\bf AMS: }
	35M33, 
	35A09, 
	35Q92, 
	65M06, 
	65M12, 
	65M20  

\section{Introduction}
Modeling and investigating the dynamics of populations is commonly viewed as one of central topics of modern mathematical demography, population biology and ecology.
Having its origin in the works of Malthus dating back to 1798 and historically preceded by Fibonacci's elementary considerations from 1202,
the mathematical theory of population dynamics underwent a rapid growth during the 19th and 20th centuries.
Among others, one should mention the works of Sharpe (1911), Lotka (1911 and 1924),  Volterra (1926), McKendrick (1926),
Kositzin (late 1930s), Fisher (1937), Kolmogorov (1937), Leslie (1945), Skellam (1950-s and 1970-s), Keyfitz (1950-s through 1980-s), 
Fredrickson \& Hoppensteadt (1971 and 1975), Gurtin (1973), Gurtin \& MacCamy (1981), etc.
For a detailed historical overview, we refer the reader to the monographs by Ianelli {\it et al.} \cite{IanMaMi2005}
and Okubo \& Levin \cite{OkuLe2001} and references therein.

The classical McKendrick-von Foerster model (often also referred to as Sharpe-Lotka-McKendrick model) reads as
\begin{equation}
	\begin{split}
		\partial_{t} p(t, a) + \partial_{a} p(t, a) &= -\mu(a) p(t, a) \text{ for } (t, a) \in (0, \infty) \times (0, a^{\dag}), \\
		p(0, a) &= p^{0}(a) \text{ for } a \in (0, a^{\dag}), \\
		p(t, 0) &= \int_{0}^{a^{\dag}} \beta(a) p(t, a) \mathrm{d}a \text{ for } a \in (0, a^{\dag})
	\end{split}
	\label{EQUATION_SHARPE_LOTKA_MCKENDRICK_MODEL}
\end{equation}
where $p(t, a)$ stands for the population individuals density of age $a \in (0, a^{\dag})$, $a^{\dag} > 0$, at time $t > 0$.
Equation (\ref{EQUATION_SHARPE_LOTKA_MCKENDRICK_MODEL}) as well as its nonlinear modifications and generalizations for the case of multiple competing populations have attracted a lot of attention.
In particular, one should mention the works and monographs by Arino \cite{Ar1992}, Chan \& Guo \cite{ChGu1989-2}, 
Ianelli {\em et al.} \cite{IanMaMi2005}, Song {\em et al.} \cite{SoYu1982}, Webb \cite{We1984}, \cite{We1985}, etc.
The questions addressed by the author range from local and global existence and uniqueness studies,
positivity and spectrum investigations as well as stability and asymptotics considerations
to optimization and control problems, etc.
The typical functional analytic framework for Equation (\ref{EQUATION_GURTIN_AND_MAC_CAMY_MODEL}) is the Lebesgue $L^{p}\big((0, a^{\dag})\big)$-space, $p \in [1, \infty)$.
Whereas most well-posedness results were obtained for $p = 1$ and similarly hold for all $p \in [1, \infty)$,
the Hilbert-space case $p = 2$ turns out to be more appropriate in some other cases (cf. \cite{Ba1991}, \cite{ChGu1989-2}).

A generalization of (\ref{EQUATION_SHARPE_LOTKA_MCKENDRICK_MODEL}) is given by Gurtin \& MacCamy's model with spatial diffusion
\begin{equation}
	\begin{split}
		\partial_{t} p(t, a, x) + \partial_{a} p(t, a, x) &= -\mu(a) p(t, a, x) + K \triangle p(t, a, x) \\
		&\phantom{=} \text{ for } (t, a, x) \in (0, \infty) \times (0, a^{\dag}) \times \Omega, \\
		p(0, a, x) &= p^{0}(a, x) \text{ for } (a, x) \in (0, a^{\dag}) \times \Omega, \\
		p(t, 0, x) &= \int_{0}^{a^{\dag}} \beta(a) p(t, a, x) \mathrm{d}a \text{ for } (t, x) \in (0, \infty) \times \Omega, \\
		p(t, a, x) &= 0 \text{ for } (t, a, x) \in (0, \infty) \times (0, a^{\dag}) \times \partial \Omega
	\end{split}
	\label{EQUATION_GURTIN_AND_MAC_CAMY_MODEL}
\end{equation}
with $p(t, a, x)$ denoting the density of the population individuals of age $a \in (0, a^{\dag})$, $a^{\dag} > 0$, at space position $x \in \Omega$
of a spatial domain $\Omega \subset \RR^{d}$ at time $t > 0$.
Global well-posedness and asymptotic behavior for Equation (\ref{EQUATION_GURTIN_AND_MAC_CAMY_MODEL}) as well as its nonlinear and stochastic versions
have been studied by Busenberg \& Iannelli \cite{Bu1983}, Chan \& Guo \cite{ChGu1989-1}, Kunisch {\em et al.} \cite{KuSchaWe1985}, Langlais \cite{La1988}, etc.
Since Equation (\ref{EQUATION_GURTIN_AND_MAC_CAMY_MODEL}) can be viewed as a ``hyperbolic-parabolic'' partial integro-differential equations,
Equation (\ref{EQUATION_GURTIN_AND_MAC_CAMY_MODEL}) is typically studied in $L^{p}\big((0, a^{\dag}) \times \Omega\big)$ for $p = 2$.

In constrast to animal populations, the migration in modern human populations is essentially nonlocal
making it possible to ignore small fluctuations arising from the random walk
and accounted for by the Laplacian term in Equation (\ref{EQUATION_GURTIN_AND_MAC_CAMY_MODEL}).
On the other hand, Equations (\ref{EQUATION_GURTIN_AND_MAC_CAMY_MODEL}) is too unrealistic
to be applied in demography since it does not account for the gender structure of the population.
To address this shortcoming, sex-structured models been developed in the 1970s, mostly within the ODE framework.
One of the first PDE models proposed is probably the one due to Keifitz.
In his article \cite[pp. 94--96]{Kei1972}, he presented a straightforward generalization of McKendrick-von Foerster model from Equation (\ref{EQUATION_SHARPE_LOTKA_MCKENDRICK_MODEL})
describing the temporal evolution of an age- and sex-structured population by the following system of partial integro-differential equations
\begin{equation}
	\begin{split}
		&\partial_{t} p_{m}(t, a_{m}) + \partial_{a_{m}} p_{m}(t, a_{m}) = -\mu_{m}(a_{m}) p_{m}(t, a_{m}) \text{ for } (t, a_{m}) \in (0, \infty) \times (0, a_{m}^{\dag}), \\
		&\partial_{t} p_{f}(t, a_{f}) + \partial_{a_{f}} p_{f}(t, a_{f}) = -\mu_{m}(a_{f}) p_{f}(t, a_{f}) \text{ for } (t, a_{f}) \in (0, \infty) \times (0, a_{f}^{\dag}), \\
		&p_{m}(0, a_{m}) = p_{m}^{0}(a_{m}) \text{ for } a_{m} \in (0, a_{m}^{\dag}), \\
		&p_{f}(0, a_{f}) = p_{f}^{0}(a_{f}) \text{ for } a_{f} \in (0, a_{f}^{\dag}), \\
		&p_{m}(0, a_{m}) =		
		\int_{0}^{a_{m}^{\dag}} \hspace{-0.2cm} \int_{0}^{a_{f}^{\dag}} \frac{\tfrac{s}{1 + s} \beta(a_{m}, a_{f})}{P_{m}(t) + P_{f}(t)} p_{m}(t, a_{m}) p_{f}(t, a_{f}) \mathrm{d}a_{m} \mathrm{d}a_{f} \text{ for }
		t \in (0, \infty) \\
		&p_{f}(0, a_{m}) =
		\int_{0}^{a_{m}^{\dag}} \hspace{-0.2cm} \int_{0}^{a_{f}^{\dag}} \frac{\tfrac{1}{1 + s} \beta(a_{m}, a_{f})}{P_{m}(t) + P_{f}(t)} p_{m}(t, a_{m}) p_{f}(t, a_{f}) \mathrm{d}a_{m} \mathrm{d}a_{f} \text{ for }
		t \in (0, \infty)
	\end{split}
	\label{EQUATION_MCKENDRICK_VON_FOERSTER_KEIFITZ_MODEL}
\end{equation}
with
\begin{equation}
	P_{m}(t) := \int_{0}^{a_{m}^{\dag}} p_{m}(t, a_{m}) \mathrm{d}a_{m}, \quad
	P_{f}(t) := \int_{0}^{a_{f}^{\dag}} p_{f}(t, a_{f}). \notag
\end{equation}
Here, $a_{m}^{\dag}, a_{f}^{\dag} \in (0, \infty]$ stand for the maximal life expectancy for male or female individuals in the population, respectively,
$p_{m}(t, a_{m})$ and $p_{f}(t, a_{f})$ denote for the number of male or female individuals of age $a_{m} \in (0, a_{m}^{\dag})$ or $a_{f} \in (0, a_{f}^{\dag})$ at time $t > 0$,
$\mu_{m}$ and $\mu_{f}$ stand for the age- and sex-specific mortality rates, $p_{m}^{0}$ and $p_{f}^{0}$ represent the population structure at the initial moment of time,
$s \in (0, 1)$ stands for the human sex ratio at birth, i.e., the ratio of male to female infants,
and $\beta(a_{m}, a_{f})$ is the birth rate in couples with a male aged $a_{m}$ and a women aged $a_{f}$.
Note that this model does not provide any information on the (official) marital status of the parents.

To account for the marital status, a new variable $c(t, a_{m}, a_{f})$
describing the number of couples with a husband of age $a_{m}$ and a wife of age $a_{f}$ at time $t$
has been introduced by Fredrickson \cite{Fre1971} and Hoppensteadt \cite{Ho1975}.
Their model is more comprehensive and contains another equation for $c$ modeling the creation and separation of couples through marriage and divorse or death
based on the so-called marriage function (see, e.g., \cite[Chapter 2.2]{IanMaMi2005}).
Whereas the necessity of incorporating the marital status into the model seemed to be very important in 1970s,
it became less significant in studying the demography of modern Western societies
due to the growing percent of single parents, childless/-free couples and singles or LGBT couples and singles
giving birth to or adapting a child.
Indeed, 40.7\% childern in the United States of America in 2011 were born to unmarried women (see \cite[p. 2]{Ma2013}) and the trend is upwards.
In 2006-2010, 43.0\% of U.S. women aged 15-44 were childless; of those who were
childless 34\% were temporarily childless,
2.3\% nonvoluntarily childless, and 6.0\% voluntarily childless (childfree) (cf. \cite[p. 4]{Ma2012}).
According to different surveys, LGBT Americans make up 3.5\%--8.0\% of the U.S. total population (see, e.g., \cite{Ga2011}).
In view of these facts, ignoring the marital status can often lead to simple and accurate demographic models.

In this article, we consider a linearized version McKendrick-von Foerster-Keifitz model from Equation (\ref{EQUATION_MCKENDRICK_VON_FOERSTER_KEIFITZ_MODEL})
which we briefly outline in Section \ref{SECTION_MODEL_DEDUCTION} below.
Then we exploit the semigroup theory to show the classical well-posedness in the sense of Hadamard
in Section \ref{SECTION_WELL_POSEDNESS_AND_ASYMPTOTICS} later on in the paper.
Under appropriate conditions on the system parameters such as fertility and mortality moduli,
we show the system to be exponentially stable.
In the subsequent Section \ref{SECTION_FINITE_DIFFERENCE_SCHEME},
we develop a finite difference scheme both with respect to age and time variables
and show it to be convergent.
Finally, in the last Section \ref{IMPLEMENATION_AND_EXAMPLE},
we discuss a computer implementation of the numerical scheme
and verify it by applying it to studying the U.S. population over the time period of 2001--2011.
Our simulation results prove to be very much consistent with the data officially reported by the U.S. Bureau of the Census \cite{USCB2013}.

\section{Model Description}
\label{SECTION_MODEL_DEDUCTION}
Let $a_{\MA}^{\dag}, a_{\VE}^{\dag} \in (0, \infty]$ be the maximal life expectancy for male or female individuals in the population, respectively.
Further, let $A_{\MA} := (0, a_{\MA}^{\dag})$, $A_{\VE} := (0, a_{\VE}^{\dag})$ be the age domains for male or female individuals, respectively.
For $t \geq 0$, let $p_{\MA}(t, a_{\MA})$ denote the total number of male individuals of age $a_{\MA} \in \bar{A}_{\MA}$ in the population.
Similarly, let $p_{\VE}(t, a_{\VE})$ denote the total number of female individuals of age $a_{\VE} \in \bar{A}_{\VE}$.
Let $\mu_{\MA}(a_{\MA})$ and $\mu_{\VE}(a_{\VE})$
be the age-specific mortality moduli of male or female individuals of age $a_{\MA} \in \bar{A}_{\MA}$ or $a_{\VE} \in \bar{A}_{\VE}$, respectively.
Further, let $b_{\MA}(a_{\MA}, a_{\VE}, p_{\MA}, p_{\VE})$ and $b_{\VE}(a_{\MA}, a_{\VE}, p_{\MA}, p_{\VE})$  describe the total number of male or female infants, respectively, born
to all couples made up of $p_{\MA}$ males of age $a_{\MA}$ and $p_{\VE}$ females of age $a_{\VE}$ with the couples being not necessarily monogamous.
Assuming
\begin{equation}
	b_{\circledast}(a_{\MA}, a_{\VE}, p_{\MA}, p_{\VE})= \int_{A_{\MA} \times A_{\VE}} \tilde{b}_{\circledast}(a_{\MA}, a_{\VE}, p_{\MA}(a_{\MA}), p_{\VE}(a_{\VE})) \mathrm{d}(a_{\MA}, a_{\VE})
	\text{ for } \circledast \in \{\MA, \VE\} \notag
\end{equation}
for some regular $\tilde{k}_{\MA}$, $\tilde{k}_{\VE}$ with
\begin{equation}
	\tilde{b}_{\MA}(\cdot, \cdot, 0, 0) \equiv \tilde{b}_{\VE}(\cdot, \cdot, 0, 0) \equiv 0 \notag
\end{equation}
and performing for each $(a_{\MA}, a_{\VE})$ a linearization of $\tilde{b}_{\circledast}(a_{\MA}, a_{\VE}, \cdot, \cdot)$ around $(0, 0)$
for $\circledast \in \{\MA, \VE\}$, we obtain the approximation
\begin{equation}
	b_{\circledast}(p_{\MA}, p_{\VE}) = \sum_{\circledcirc \in \{\MA, \VE\}}
	\int_{A_{\circledcirc}} \beta_{\circledast \circledcirc}(a_{\circledcirc}) \mathrm{d}a_{\circledcirc}
	\text{ for } \circledast \in \{\MA, \VE\} \notag
\end{equation}
with
\begin{align*}
	\beta_{\MA \MA}(a_{\MA}) &= \int_{A_{\VE}} \partial_{p_{\MA}} \tilde{b}_{\MA}(a_{\MA}, a_{\VE}, 0, 0) \mathrm{d}a_{\VE}, & \beta_{\MA \VE}(a_{\VE}) &= \int_{A_{\MA}} \partial_{b_{\VE}} \tilde{b}_{\MA}(a_{\MA}, a_{\VE}, 0, 0) \mathrm{d}a_{\MA}, \\
	\beta_{\VE \MA}(a_{\MA}) &= \int_{A_{\VE}} \partial_{p_{\MA}} \tilde{b}_{\VE}(a_{\MA}, a_{\VE}, 0, 0) \mathrm{d}a_{\VE}, & \beta_{\VE \VE}(a_{\VE}) &= \int_{A_{\MA}} \partial_{b_{\VE}} \tilde{b}_{\MA}(a_{\MA}, a_{\VE}, 0, 0) \mathrm{d}a_{\MA}.
\end{align*}
Here, $\beta_{\MA \MA}(a_{\MA}), \dots, \beta_{\VE \VE}(a_{\VE})$ stand for the age- and sex-specific fertility moduli for male or female infants.
Usually, $\beta_{\MA \MA} \equiv \beta_{\VE \MA} \equiv 0$, $0 < \beta_{\MA \VE} \approx \beta_{\VE \VE} > 0$
since the influence of the male part of population is overwhelmingly nonlinear (cf. \cite{RaKo2007}).
Further, let $g_{\MA}(t, a_{\MA})$ and $g_{\VE}(t, a_{\VE})$
be the net immigration of male or female individuals of age $a_{\MA} \in A_{\MA}$ or $a_{\VE} \in A_{\VE}$, respectively, at time $t > 0$.
With $p_{\MA}^{0}(a_{\MA})$ and $p_{\VE}^{0}(a_{\VE})$ describing the total number of male or female individuals of age $a_{\MA} \in \bar{A}_{\MA}$ or $a_{\VE} \in \bar{A}_{\VE}$, respectively, 
in the population at the initial moment of time and
$g_{\MA}(t, a_{\MA})$ and $g_{\VE}(t, a_{\VE})$ quantifying the net immigration of male or female individuals of age $a_{\MA} \in \bar{A}_{\MA}$ or $a_{\VE} \in \bar{A}_{\VE}$ at time $t$, 
the evolution equations for $(p_{\MA}, p_{\VE})$ read as
\begin{align}
	\partial_{t} p_{\MA}(t, a_{\MA}) + \partial_{a_{\MA}} p_{\MA}(t, a_{\MA}) &= -\mu_{\MA}(a_{\MA}) p_{\MA}(t, a_{\MA}) + g_{\MA}(t, a_{\MA}) \notag \\
	&\phantom{=} \text{ for } (t, a_{\MA}) \in (0, \infty) \times A_{\MA},
	\label{EQUATION_MODEL_EQUATION_1} \\
	\partial_{t} p_{\VE}(t, a_{\VE}) + \partial_{a_{\VE}} p_{\VE}(t, a_{\VE}) &= -\mu_{\VE}(a_{\VE}) p_{\VE}(t, a_{\VE}) + g_{\VE}(t, a_{\VE}) \notag \\
	&\phantom{=} \text{ for } (t, a_{\MA}) \in (0, \infty) \times A_{\MA},
	\label{EQUATION_MODEL_EQUATION_2} \\
	p_{\MA}(t, 0) &= \sum_{\circledcirc \in \{\MA, \VE\}} \int_{A_{\circledcirc}} \beta_{\MA \circledcirc}(a_{\circledcirc}) p_{\circledcirc}(t, a_{\circledcirc}) \mathrm{d}a_{\circledcirc} \text{ for } t \in (0, \infty),
	\label{EQUATION_MODEL_EQUATION_3} \\
	p_{\VE}(t, 0) &= \sum_{\circledcirc \in \{\MA, \VE\}} \int_{A_{\circledcirc}} \beta_{\VE \circledcirc}(a_{\circledcirc}) p_{\circledcirc}(t, a_{\circledcirc}) \mathrm{d}a_{\circledcirc} \text{ for } t \in (0, \infty),
	\label{EQUATION_MODEL_EQUATION_4} \\
	p_{\MA}(0, a_{\MA}) &= p_{\MA}^{0}(a_{\MA}) \text{ for } a_{\MA} \in A_{\MA},
	\label{EQUATION_MODEL_EQUATION_5} \\
	p_{\VE}(0, a_{\VE}) &= p_{\VE}^{0}(a_{\VE}) \text{ for } a_{\VE} \in A_{\VE}.
	\label{EQUATION_MODEL_EQUATION_6}
\end{align}
Here, Equations (\ref{EQUATION_MODEL_EQUATION_1})--(\ref{EQUATION_MODEL_EQUATION_2})
represent a conservation law describing the natural ageing and migration
whereas Equations (\ref{EQUATION_MODEL_EQUATION_3})--(\ref{EQUATION_MODEL_EQUATION_4})
stand for the so-called ``birth law''
being a boundary condition with a non-local term.
Finally, Equations (\ref{EQUATION_MODEL_EQUATION_5})--(\ref{EQUATION_MODEL_EQUATION_6}) prescribe the initial population structure.

Following \cite{We1985}, we assume $\mu_{\MA} \colon A_{\MA} \to [0, \infty), \mu_{\VE} \colon A_{\VE} \to [0, \infty)$ to be Lebesgue-integrable and
define the survival probability for male or female individuals till the age $a_{\MA} \in \bar{A}_{\MA}$ or $a_{\VE} \in \bar{A}_{\VE}$, respectively, as
\begin{equation}
	\pi_{\MA}(a_{\MA}) := \exp\Big(-\int_{0}^{a_{\MA}} \mu_{\MA}(a_{\MA}) \mathrm{d}a_{\MA}\Big) \text{ and }
	\pi_{\VE}(a_{\VE}) := \exp\Big(-\int_{0}^{a_{\VE}} \mu_{\VE}(a_{\VE}) \mathrm{d}a_{\VE}\Big). \notag
\end{equation}
For $\pi_{\MA}$ and $\pi_{\VE}$ to vanish in $a_{\MA}^{\dag}$ or $a_{\VE}^{\dag}$, respectively,
we require that the integrals
\begin{equation}
	\int_{A_{\MA}} \mu_{\MA}(a_{\MA}) \mathrm{d}a_{\MA} = \infty \text{ and }
	\int_{A_{\VE}} \mu_{\VE}(a_{\VE}) \mathrm{d}a_{\VE} = \infty \notag
\end{equation}
are divergent. For finite $a_{\MA}^{\dag}, a_{\VE}^{\dag}$, this would mean $\mu_{\MA} \not \in L^{\infty}(A_{\MA})$, $\mu_{\VE} \not \in L^{\infty}(A_{\VE})$.
In contrast to that, we have $\pi_{\MA} \in L^{\infty}(A_{\MA})$, $\pi_{\VE} \in L^{\infty}(A_{\VE})$
both for finite and infnite $a_{\MA}^{\dag}, a_{\VE}^{\dag}$.
Additionally, we impose the natural condition
\begin{equation}
	m_{\MA \MA} \in L^{\infty}(A_{\MA}), \dots, m_{\VE \VE} \in L^{\infty}(A_{\VE}). \notag
\end{equation}
The latter is satisfied if $\beta_{\MA \MA} \in L^{\infty}(A_{\MA}), \dots, \beta_{\VE \VE} \in L^{\infty}(A_{\VE})$
exhibit a sufficiently rapid decay rate in $a_{\MA}^{\dag}$ or $a_{\VE}^{\dag}$, respectively.

Thus, to avoid the necessity of working with weighted Lebesgue- and Sobolev spaces,
similar to \cite[p. 255]{CuGeGio2009}, we define the new variables
\begin{equation}
	u_{\MA}(t, a_{\MA}) := \tfrac{p_{\MA}(t, a_{\MA})}{\pi_{\MA}(a_{\MA})} \text{ and }
	u_{\VE}(t, a_{\VE}) := \tfrac{p_{\VE}(t, a_{\VE})}{\pi_{\VE}(a_{\VE})} \text{ for } a_{\MA} \in A_{\MA}, a_{\VE} \in A_{\VE}. \notag
\end{equation}
Introducing the age- and sex-specific maternity functions
\begin{equation}
	\begin{split}
		m_{\MA \MA}(a_{\MA}) &:= \pi_{\MA}(a_{\MA}) \beta_{\MA \MA}(a_{\MA}), \quad
		m_{\MA \VE}(a_{\VE}) := \pi_{\MA}(a_{\MA}) \beta_{\MA \VE}(a_{\VE}), \\
		m_{\VE \MA}(a_{\MA}) &:= \pi_{\VE}(a_{\VE}) \beta_{\VE \MA}(a_{\MA}), \quad
		m_{\VE \VE}(a_{\VE}) := \pi_{\VE}(a_{\VE}) \beta_{\VE \VE}(a_{\VE}),
	\end{split}
	\notag
\end{equation}
we can use Equations (\ref{EQUATION_MODEL_EQUATION_1})--(\ref{EQUATION_MODEL_EQUATION_6}) to easily verify that $(u_{\MA}, u_{\VE})$ solves the problem
\begin{align}
	\partial_{t} u_{\MA}(t, a_{\MA}) + \partial_{a_{\MA}} u_{\MA}(t, a_{\MA}) &= f_{\MA}(t, a_{\MA}) \text{ for } (t, a_{\MA}) \in (0, \infty) \times A_{\MA},
	\label{EQUATION_MODEL_TRANSFORMED_EQUATION_1} \\
	\partial_{t} p_{\VE}(t, a_{\VE}) + \partial_{a_{\VE}} u_{\VE}(t, a_{\VE}) &= f_{\VE}(t, a_{\VE}) \text{ for } (t, a_{\VE}) \in (0, \infty) \times A_{\VE},
	\label{EQUATION_MODEL_TRANSFORMED_EQUATION_2} \\
	u_{\MA}(t, 0) &= \sum_{\circledcirc \in \{\MA, \VE\}} \int_{A_{\circledcirc}} m_{\MA \circledcirc}(a_{\circledcirc}) u_{\circledcirc}(t, a_{\circledcirc}) \mathrm{d}a_{\circledcirc} \text{ for } t \in (0, \infty),
	\label{EQUATION_MODEL_TRANSFORMED_EQUATION_3} \\
	u_{\VE}(t, 0) &= \sum_{\circledcirc \in \{\MA, \VE\}} \int_{A_{\circledcirc}} m_{\VE \circledcirc}(a_{\circledcirc}) u_{\circledcirc}(t, a_{\circledcirc}) \mathrm{d}a_{\circledcirc} \text{ for } t \in (0, \infty),
	\label{EQUATION_MODEL_TRANSFORMED_EQUATION_4} \\
	u_{\MA}(0, a_{\MA}) &= u_{\MA}^{0}(a_{\MA}) \text{ for } a_{\MA} \in A_{\MA},
	\label{EQUATION_MODEL_TRANSFORMED_EQUATION_5} \\
	u_{\VE}(0, a_{\VE}) &= u_{\VE}^{0}(a_{\VE}) \text{ for } a_{\VE} \in A_{\VE},
	\label{EQUATION_MODEL_TRANSFORMED_EQUATION_6}
\end{align}
where
\begin{equation}
    u_{\MA}^{0}(a_{\MA}) := \tfrac{p_{\MA}^{0}(a_{\MA})}{\pi_{\MA}(a_{\MA})}, \quad
	u_{\VE}^{0}(a_{\VE}) := \tfrac{p_{\VE}^{0}(a_{\VE})}{\pi_{\VE}(a_{\VE})}
    \text{ for } a_{\MA} \in A_{\MA}, a_{\VE} \in A_{\VE} \notag
\end{equation}
and
\begin{equation}
    f_{\MA}(t, a_{\MA}) := \tfrac{g_{\MA}(t, a_{\VE})}{\pi_{\MA}(a_{\MA})}, \quad
    f_{\VE}(t, a_{\VE}) := \tfrac{g_{\VE}(t, a_{\VE})}{\pi_{\VE}(a_{\VE})}
    \text{ for } t > 0, a_{\MA} \in A_{\MA}, a_{\VE} \in A_{\VE}. \notag
\end{equation}

\section{Well-posedness and Long-Time Behavior}
\label{SECTION_WELL_POSEDNESS_AND_ASYMPTOTICS}
In this section, we want to prove the classical well-posedness in the sense of Hadamard for (\ref{EQUATION_MODEL_TRANSFORMED_EQUATION_1})--(\ref{EQUATION_MODEL_TRANSFORMED_EQUATION_6}).
To this end, we state the problem in a Hilbert space setting and apply the operator semigroup theory (see \cite{ArBaHieNeu2001}, \cite{Pa1983}).
Our approach differs inasmuch from the classical one (see, e.g., \cite{We1985} and references therein)
as we use the semigroup theory instead of Fredholm integral equation theory to obtain the well-posedness.
Further, unlike other authors (cf. \cite{Ar1992}, \cite{We1984}) who also applied the semigroup theory to similar problems,
we exploit only Hilbert space techniques rather then working with the $L^{1}$-space.
Though at first glance the $L^{2}$-space might appear to be not the most intuitive choice since 
it the $L^{2}$-norm can not be directly related to the population size,
it provides more structure and thus facilitates the analytical and numerical treatment of the problem
without being an actual restriction in demographical applications.

In the following, we assume $m_{\MA \MA}, m_{\VE \MA} \in L^{\infty}(A_{\MA})$ and $m_{\MA \VE}, m_{\VE \VE} \in L^{\infty}(A_{\VE})$.
We consider the Hilbert space $X := L^{2}(A_{\MA}) \times L^{2}(A_{\VE})$ endowed with the standard product topology.
We define the operator $A \colon D(A) \subset X \to X$ given as
\begin{equation}
	(u_{\MA}(a_{\MA}), u_{\VE}(u_{\VE})) =: u \mapsto
	\begin{pmatrix}
		-\partial_{a_{\MA}} u_{\MA}(a_{\MA}) \\
		-\partial_{a_{\VE}} u_{\VE}(a_{\VE})
	\end{pmatrix}
	\text{ for } a_{\MA} \in A_{\MA}, a_{\VE} \in A_{\VE}
	\notag
\end{equation}
with the domain
\begin{equation}
	\begin{split}
		D(A) := \Big\{(u_{\MA}, u_{\VE}) \in H^{1}(A_{\MA}) \times H^{1}(A_{\VE}) \,\Big|\,
		&u_{\MA}(0) = \sum_{\circledcirc \in \{\MA, \VE\}} \int_{A_{\circledcirc}} m_{\MA \circledcirc}(a_{\circledcirc}) u_{\circledcirc}(a_{\circledcirc}) \mathrm{d}a_{\circledcirc}, \\
		&u_{\VE}(a) = \sum_{\circledcirc \in \{\MA, \VE\}} \int_{A_{\circledcirc}} m_{\VE \circledcirc}(a_{\circledcirc}) u_{\circledcirc}(a_{\circledcirc}) \mathrm{d}a_{\circledcirc}\Big\}
	\end{split}
	\notag
\end{equation}
equipped with the standard product topology on $H^{1}(A_{\MA}) \times H^{1}(A_{\VE})$.
Here and in the sequel, $H^{1} := H^{1, 2} = W^{1, 2}$ will denote the standard scalar-valued (see, e.g., \cite[Chapter 3]{Ad1975})
or Banach-space-valued Sobolev space (cf, e.g., \cite[p. 2]{Si1990}).

\begin{remark}
	Under the condition
	\begin{equation}
		\begin{split}
			\Big(1 - \int_{A_{\MA}} m_{\MA \MA}(a_{\MA}) \mathrm{d}a_{\MA}\Big)
			&\Big(1 - \int_{A_{\VE}} m_{\VE \VE}(a_{\VE}) \mathrm{d}a_{\VE}\Big) \neq \\
			&\Big(\int_{A_{\MA}} m_{\MA \VE}(a_{\VE}) \mathrm{d}a_{\VE}\Big)
			\Big(\int_{A_{\VE}} m_{\VE \MA}(a_{\MA}) \mathrm{d}a_{\MA}\Big),
		\end{split}
		\notag
	\end{equation}
	the expression
	\begin{equation}
		\Big[\sum_{\circledast \in \{\MA, \VE\}} \Big(u_{\circledast}(0) - \sum_{\circledcirc \in \{\MA, \VE\}} \int_{A_{\circledcirc}} m_{\circledast\circledcirc}(a_{\circledcirc})
		u_{\circledcirc}(a_{\circledcirc}) \mathrm{d}a_{\circledcirc}\Big)^{2}\Big]^{1/2}
		\notag
	\end{equation}
	gives a seminorm on $H^{1}(A_{\MA}) \times H^{1}(A_{\VE})$, being additionally a norm on the subspace of constant functions, and thus
	\begin{equation}
		u \mapsto \Big(\sum_{\circledast \in \{\MA, \VE\}} \|\partial_{a_{\circledast}} u_{\circledast}\|_{L^{2}(A_{\circledast})}^{2}\Big)^{1/2} =
		\|A u\|_{X} \notag
	\end{equation}
	constitutes an equivalent norm on $D(A)$ by virtue of the third Poincar\'{e}'s inequality.
\end{remark}

Due to the Sobolev embedding theory (cf. \cite[Theorem 4.12]{Ad1975}), we know
\begin{equation}
	H^{1}(A_{\MA}) \hookrightarrow C^{0}_{b}(\bar{A}_{\MA}), \quad
	H^{1}(A_{\VE}) \hookrightarrow C^{0}_{b}(\bar{A}_{\VE}). \notag
\end{equation}
Thus, $A$ is well-defined. The linearity of $A$ is also obvious.

Letting, $u := (u_{\MA}, u_{\VE})$, $f := (f_{\MA}, f_{\VE})$
and $u^{0} = (u_{\MA}^{0}, u_{\VE}^{0})$, Equations (\ref{EQUATION_MODEL_TRANSFORMED_EQUATION_1})--(\ref{EQUATION_MODEL_TRANSFORMED_EQUATION_6}) can now be equivalently written in the abstract form
\begin{equation}
	\dot{u}(t) = A u(t) + f(t) \text{ for } t > 0, \quad u(0) = u^{0}.
	\label{EQUATION_MODEL_ABSTRACT_FORM}
\end{equation}

Since we will observe that $A$ is closed and has a non-empty resolvent set 
(cf. Lemmas \ref{LEMMA_OPERATOR_A_DENSE_AND_CLOSED} and \ref{LEMMA_A_MINUS_BETA_M_DISSIPATIVE} below),
by a well-known result on operator semigroups (see, e.g., \cite[Theorem 3.1.12]{ArBaHieNeu2001}),
proving the classical well-posedness for the abstract Cauchy problem (\ref{EQUATION_MODEL_ABSTRACT_FORM})
and thus also for the original initial-boundary value problem (\ref{EQUATION_MODEL_EQUATION_1})--(\ref{EQUATION_MODEL_EQUATION_6})
reduces to showing that $A$ is an infinitesimal generator of $C_{0}$-semigroup of bounded linear operators on $X$.

\begin{lemma}
	\label{LEMMA_OPERATOR_A_DENSE_AND_CLOSED}
	The operator $A$ is densely defined and closed.
\end{lemma}

\begin{proof}
	\begin{itemize}
		\item[] \emph{Density:} Let $u := (u_{\MA}, u_{\VE}) \in X$ and let $\varepsilon > 0$ be arbitrary.
		Due to the density of test functions in $L^{2}(A_{\MA})$ and $L^{2}(A_{\VE})$ as well as the monotonicity of Lebesgue integral,
		there exists a number $\delta_{0} \in \big(0, \max\{a_{\MA}^{\dag}, a_{\VE}^{\dag}\}\big)$ such that for any $\delta \in (0, \delta_{0})$ there exist test functions
		$\varphi_{\MA}(\cdot; \delta) \in C_{0}^{\infty}(A_{\MA})$, $\varphi_{\VE}(\cdot; \delta) \in C_{0}^{\infty}(A_{\VE})$ such that
		\begin{equation}
			\mathrm{supp}\,(\varphi_{\circledast}(\cdot; \delta)) \subset (\delta, a_{\circledast}^{\dag}) \text{ and }
			\|u_{\circledast} - \varphi_{\circledast}(\cdot; \delta)\|_{L^{2}(A_{\circledast})} < \varepsilon/2 \text{ for } \circledast \in \{\MA, \VE\}. \notag
		\end{equation}
		We let
		\begin{equation}
			I_{\circledast}(\delta) := \sum_{\circledcirc \in \{\MA, \VE\}} \int_{A_{\circledcirc}} m_{\circledast \circledcirc}(a_{\circledcirc}) \varphi_{\circledcirc}(a_{\circledcirc}; \delta) \mathrm{d}a_{\circledcirc}
			\text{ for } \circledast \in \{\MA, \VE\}.
			\notag
		\end{equation}
		Note that by the virtue of H\"older's inequality,
		both $I_{\MA}(\delta)$ and $I_{\VE}(\delta)$ are absolutely and uniformly bounded with respect to $\delta \in (0, \delta_{0})$ by the number
		\begin{equation}
			C := 2 M \max\big\{(a_{\MA}^{\dag})^{1/2}, (a_{\VE}^{\dag})^{1/2}\big\} \max\{\|u_{\MA}\|_{L^{2}(A_{\MA})} + \tfrac{\varepsilon}{2},
			\|u_{\VE}\|_{L^{2}(A_{\VE})} + \tfrac{\varepsilon}{2}\big\}
			\notag
		\end{equation}
		with
		\begin{equation}
			M := \max_{\circledast, \circledcirc \in \{\MA, \VE\}}  \|m_{\circledast \circledcirc}\|_{L^{\infty}(A_{\circledcirc})}. \notag
		\end{equation}

		For $a \geq 0$, $\delta \in (0, \delta_{0})$ and $\theta \in \RR$, consider the measurable function
		\begin{equation}
			l(a; \delta, \theta) := \tfrac{\theta}{\delta}(\delta - a) \chi_{[0, \; \delta]}(a) \notag
		\end{equation}
		with $\chi_{[0, \delta]}$ standing for the characteristic function of $[0, \delta]$. Letting
		\begin{equation}
			u_{\circledast, \varepsilon}(a_{\circledast}; \delta, \theta_{\circledast}) := \varphi_{\circledast}(a_{\circledast}; \delta) +
			l(a_{\circledast}; \delta, \theta_{\circledast}) \text{ for } a_{\circledast} \in A_{\circledast}, \circledast \in \{\MA, \VE\},
			\notag
		\end{equation}
		we observe $u_{\MA, \varepsilon} \in H^{1}(A_{\MA})$, $u_{\VE, \varepsilon} \in H^{1}(A_{\VE})$.
		Now, the parameters $\delta \in (0, \delta_{0})$, $\theta_{\MA}, \theta_{\VE} \in \RR$ have to be selected such that
		\begin{equation}
			\begin{split}
				\big(u_{\MA, \varepsilon}(\cdot; \delta, \theta_{\MA}), u_{\VE, \varepsilon}(\cdot; \delta, \theta_{\VE})\big) &\in D(A) \text{ and } \\
				\Big\|\big(u_{\MA, \varepsilon}(\cdot; \delta, \theta_{\MA}), u_{\VE, \varepsilon}(\cdot; \delta, \theta_{\VE})\big) - (u_{\MA}, u_{\VE})\Big\|_{X} &< \varepsilon
			\end{split}
			\notag
		\end{equation}
		holds true, i.e., there suffices to fulfil
		\begin{equation}
			\begin{split}
				u_{\circledast, \varepsilon}(0; \delta, \theta_{\circledast}) = \sum_{\circledcirc \in \{\MA, \VE\}}
				\int_{A_{\circledcirc}} m_{\circledast \circledcirc}(a_{\circledcirc}) \varphi_{\circledcirc}(a_{\circledcirc}; \delta) \mathrm{d}a_{\circledcirc}
				\text{ for } \circledast \in \{\MA, \VE\} \text{ and } \\
				\int_{A_{\circledast}} |u_{\circledast, \varepsilon}(a_{\circledast}; \delta, \theta_{\circledast}) - u_{\circledast}(a_{\circledast})|^{2} \mathrm{d}a_{\circledast} < \varepsilon^{2}/2
				\text{ for } \circledast \in \{\MA, \VE\}.
			\end{split}
			\notag
		\end{equation}
		The latter conditions are satisfied if
		\begin{align}
			\theta_{\MA} =
			\sum_{\circledast \in \{\MA, \VE\}} \int_{A_{\circledast}} m_{\MA \circledast}(a_{\circledast}) l(a_{\circledast}; \delta, \theta_{\circledast}) \mathrm{d}a_{\circledast} + I_{\MA}(\delta),
			\label{EQUATION_CONDITION_DELTA_THETA_1} \\
			|\theta_{\MA}| \sqrt{\delta}/\sqrt{3} < \varepsilon^{2}/4,
			\label{EQUATION_CONDITION_DELTA_THETA_2} \\
			\theta_{\VE} =
			\sum_{\circledast \in \{\MA, \VE\}} \int_{A_{\circledast}} m_{\VE \circledast}(a_{\circledast}) l(a_{\circledast}; \delta, \theta_{\circledast}) \mathrm{d}a_{\circledast} + I_{\VE}(\delta),
			\label{EQUATION_CONDITION_DELTA_THETA_3} \\
			|\theta_{\VE}| \sqrt{\delta}/\sqrt{3} < \varepsilon^{2}/4.
			\label{EQUATION_CONDITION_DELTA_THETA_4}
		\end{align}
		Estimating for $\circledast \in \{\MA, \VE\}$
		\begin{equation}
			\begin{split}
				\Big|\sum_{\circledcirc \in \{\MA, \VE\}} \int_{A_{\circledcirc}} m_{\circledast \circledcirc}(a_{\circledcirc}) l(a_{\circledcirc}; \delta, \theta_{\circledcirc}) \mathrm{d}a_{\circledcirc}\Big| &\leq
				M \delta \max\{\theta_{\MA}, \theta_{\VE}\} \max\big\{(a_{\MA}^{\dag})^{1/2}, (a_{\VE}^{\dag})^{1/2}\big\} \\
				&=: \kappa(\delta) \max\{\theta_{\MA}, \theta_{\VE}\}
			\end{split}
			\notag
		\end{equation}
		and observing that the matrix
		$\begin{pmatrix}
			1 + \kappa_{11}(\delta) & \kappa_{12}(\delta) \\
			\kappa_{21}(\delta) & 1 + \kappa_{22}(\delta)
		\end{pmatrix}$
		is invertible with the operator norm of the inverse being uniformly bounded by $3$ 
		if $\max\limits^{}_{i, j = 1, 2} |\kappa_{ij}(\delta)| \leq \tfrac{1}{3}$, i.e., if, e.g.,
		\begin{equation}
			\delta \in (0, \min\{\delta_{0}, \delta_{1}\}) \text{ with }
			\delta_{1} := \tfrac{3}{(M + 1) \max\{\theta_{\MA}, \theta_{\VE}\} \max\big\{(a_{\MA}^{\dag})^{1/2}, (a_{\VE}^{\dag})^{1/2}\big\}}, \notag
		\end{equation}
		we conclude that the linear system (\ref{EQUATION_CONDITION_DELTA_THETA_1}), (\ref{EQUATION_CONDITION_DELTA_THETA_3})
		is uniquely solvable for $(\theta_{\MA}, \theta_{\VE})$ with
		\begin{equation}
			\theta_{\MA}^{2} + \theta_{\VE}^{2} \leq 9 \big(I_{\MA}^{2}(\delta) + I_{\VE}^{2}(\delta)\big) \leq 18 C^{2}. \notag
		\end{equation}
		Hence, selecting
		\begin{equation}
			\delta \in \min\{\delta_{0}, \delta_{1}, \delta_{2}\} \text{ with } \delta_{2} := \tfrac{\varepsilon^{4}}{24(1 + C^{2})}, \notag
		\end{equation}
		all equations and inequalities in (\ref{EQUATION_CONDITION_DELTA_THETA_1})--(\ref{EQUATION_CONDITION_DELTA_THETA_4})
		are satisfied.
		Thus, the constructed function $\big(u_{\MA, \varepsilon}(\cdot; \delta, \theta_{\MA}), u_{\VE, \varepsilon}(\cdot; \delta, \theta_{\VE})\big) \in D(A)$
		lies in an $\varepsilon$-neighborhood of $u$.

		\item[] \emph{Closedness:}
		We consider the operator $F \colon H^{1}(A_{\MA}) \times H^{1}(A_{\VE}) \to \RR \times \RR$ with
		\begin{equation}
			(u_{\MA}, u_{\VE}) \mapsto
			\begin{pmatrix}
				u_{\MA}(0) - \sum\limits_{\circledast \in \{\MA, \VE\}} \int_{A_{\circledast}} m_{\MA \circledast}(a_{\circledast}) \varphi_{\circledast}(a_{\circledast}; \delta) \mathrm{d}a_{\circledast} \\
				u_{\VE}(0) - \sum\limits_{\circledast \in \{\MA, \VE\}} \int_{A_{\circledast}} m_{\VE \circledast}(a_{\circledast}) \varphi_{\circledast}(a_{\circledast}; \delta) \mathrm{d}a_{\circledast}
			\end{pmatrix}.
			\notag
		\end{equation}
		By the virtue of Sobolev embedding theorem, $F$ is a bounded linear operator.
		Since $\{(0, 0)\}$ is a closed subspace of $\RR \times \RR$
		and $D(A) = F^{-1}\big(\{(0, 0)\}\big)$, the latter is a closed subspace of $H^{1}(A_{\MA}) \times H^{1}(A_{\VE})$
		and thus a Banach space.
		Now, the operator $A$ is bounded linear map between the Banach spaces $D(A)$ and $X$ and therefore a closed linear operator.
	\end{itemize}
	The proof is finished.
\end{proof}

\begin{lemma}
	\label{LEMMA_A_MINUS_BETA_M_DISSIPATIVE}
	For $\beta > 0$ sufficiently large, the operator $A - \beta \mathrm{id}_{X}$ is m-dissipative.
\end{lemma}

\begin{proof}
	For $\beta \in \mathbb{R}$ and $u \in D(A)$, we have
	\begin{equation}
		\begin{split}
			\langle (A - \beta \mathrm{id}_{X}) u, u\rangle_{X} &=
			\langle A u, u\rangle_{X} - \beta \|u\|_{X}^{2} \\
			&= -\langle \partial_{a_{\MA}} u_{\MA}, u_{\MA}\rangle_{L^{2}(A_{\MA})}
			-\langle \partial_{a_{\VE}} u_{\VE}, u_{\VE}\rangle_{L^{2}(A_{\VE})} - \beta \|u\|_{X}^{2} \\
			&\leq -u_{\MA}^{2}(a_{\VE}^{\dag}) + u_{\MA}^{2}(0) - u_{\VE}^{2}(a_{\VE}^{\dag}) + u_{\VE}^{2}(0) - \beta \|u\|_{X}^{2} \\
			&\leq
			-\sum_{\circledast \in \{\MA, \VE\}}
			\Big(\sum_{\circledcirc \in \{\MA, \VE\}} \int_{A_{\circledcirc}} m_{\circledast \circledcirc}(a_{\circledcirc}) u_{\circledcirc}(a_{\circledcirc}) \mathrm{d}a_{\circledcirc}\Big)^{2} - \beta \|u\|_{X}^{2} \\
			&\leq -\Big(\beta - \sum_{\circledast, \circledcirc \in \{\MA, \VE\}} 2 a_{\circledcirc}^{\dag} \|m_{\circledast \circledcirc}\|_{L^{\infty}(A_{\circledcirc})}^{2}\Big) \|u\|_{X}^{2}.
		\end{split}
		\label{EQUATION_DISSIPATIVITY_OF_A_MINUS_BETA}
	\end{equation}
	Thus, for $\beta \geq \beta_{0}$ with
	\begin{equation}
		\beta_{0} := \sum_{\circledast, \circledcirc \in \{\MA, \VE\}} 2 a_{\circledcirc}^{\dag} \|m_{\circledast \circledcirc}\|_{L^{\infty}(A_{\circledcirc})}^{2}, \notag
	\end{equation}
	the operator $A - \beta \mathrm{id}_{X}$ is dissipative.
	
	Next, we show that the operator $A - \beta \mathrm{id}_{X}$ is surjective for some $\beta \geq \beta_{0}$.
	For $f = (f_{\MA}, f_{\VE})$, we solve for $u \in D(A)$ the equation
	\begin{equation}
		(A - \beta \mathrm{id}_{X}) u = f.
		\label{EQUATION_OPERATOR_EQUATION}
	\end{equation}
	Multiplying Equation (\ref{EQUATION_OPERATOR_EQUATION}) with $v$ in $X$, we obtain the weak formulation
	\begin{equation}
		a(u, v; \beta) = -\langle f, v\rangle_{X} \text{ for all } v \in X \label{EQUATION_ELLIPTIC_EQUATION_WEAK_FORMULATION}
	\end{equation}
	with the bilinear form $a(\cdot, \cdot; \beta) \colon D(A) \times X \to \RR$ given as
	\begin{equation}
		a(u, v; \beta) := \sum_{\circledast \in \{\MA, \VE\}} \int_{A_{\circledast}}
		\bigg(\Big(\partial_{a_{\circledast}} u(a_{\circledast}) v(a_{\circledast})\big) + \beta u(a_{\circledast}) v(a_{\circledast})\Big) \mathrm{d}a_{\circledast}
		\text{ for } u \in D(A), v \in X.
		\notag
	\end{equation}
	Now, we want to apply Babu\v{s}ka-Lax-Milgram lemma to solve Equation (\ref{EQUATION_ELLIPTIC_EQUATION_WEAK_FORMULATION}).
	This amounts to showing that $a(\cdot, \cdot; \beta)$ is continuous on $D(A) \times X$
	and satisfies the inf-sup condition
	\begin{equation}
		\inf_{u \in D(A)} \sup_{v \in X} \tfrac{a(u, v; \beta)}{\|u\|_{D(A)} \|v\|_{X}} > 0. \notag
	\end{equation}
	Whereas the continuity of $a(\cdot, \cdot; \beta)$ is obvious, the inf-sup-condition
	holds true if and only if there exist constants $c_{1}(\beta), c_{2}(\beta) > 0$ such that
	for any $v \in X$ there exists $u \in D(A)$ such that
	\begin{equation}
		a(u, v; \beta) \geq c_{1}(\beta) \|v\|_{X}^{2} \text{ and } \|u\|_{D(A)} \leq c_{2}(\beta) \|v\|_{X}. \notag
	\end{equation}
	Indeed, let $v \in V$ be arbitrary. For a sufficiently large $\beta$, we look for $u \in D(A)$ satisfying
	\begin{equation}
		\partial_{a_{\circledast}} u_{\circledast} + \beta u_{\circledast} = v_{\circledast} \text{ for } \circledast \in \{\MA, \VE\}
		\label{EQUATION_CONSTRUCTION_OF_U_FOR_BABUSKA_LAX_MILGRAM}
	\end{equation}
	where the condition $u \in D(A)$ dictates
	\begin{equation}
		u_{\circledast}(0) = \sum_{\circledcirc \in \{\MA, \VE\}} \int_{A_{\circledcirc}} m_{\circledast \circledcirc}(a_{\circledcirc}) u_{\circledcirc}(a_{\circledcirc}) \mathrm{d}a_{\circledcirc} \text{ for } \circledast \in \{\MA, \VE\}.
		\label{EQUATION_CONSTRUCTION_OF_U_FOR_BABUSKA_LAX_MILGRAM_IC}
	\end{equation}
	From Equation (\ref{EQUATION_CONSTRUCTION_OF_U_FOR_BABUSKA_LAX_MILGRAM}), we obtain by the virtue of Duhamel's formula
	\begin{equation}
		u_{\circledast}(a_{\circledast}) = c_{\circledast} e^{-\beta a_{\circledast}} +
		\int_{0}^{a_{\circledast}} e^{-\beta(a_{\circledast} - \alpha_{\circledast})} v_{\circledast}(\alpha_{\circledast}) \mathrm{d}\alpha_{\circledast}
		\label{EQUATION_CONSTRUCTION_OF_U_FOR_BABUSKA_LAX_MILGRAM_SOLUTION}
	\end{equation}
	for some constants $c_{\MA}, c_{\VE} \in \RR$.
	Note that we trivially have $u_{\circledast} \in H^{1}(A_{\circledast})$ since
	\begin{equation}
		\|u\|_{H^{1}(A_{\circledast}) \times H^{1}(A_{\circledcirc})}^{2} \leq
		2 (1 + \beta) \max\{a_{\MA}^{\dag}, a_{\VE}^{\dag}\} \big(\|c\|_{\RR^{2}}^{2} + \|v\|_{X}^{2}\big) =: C_{1}^{2}(\beta) \big(\|c\|_{\RR^{2}}^{2} + \|v\|_{X}^{2}\big). \notag
	\end{equation}
	Equations (\ref{EQUATION_CONSTRUCTION_OF_U_FOR_BABUSKA_LAX_MILGRAM_SOLUTION}), (\ref{EQUATION_CONSTRUCTION_OF_U_FOR_BABUSKA_LAX_MILGRAM_IC}) yield a linear system for $(c_{\MA}, c_{\VE})$
	\begin{equation}
		c_{\circledast} = \sum_{\circledcirc \in \{\MA, \VE\}} \int_{A_{\circledcirc}} m_{\circledast \circledcirc}(a_{\circledcirc}) e^{-\beta a_{\circledcirc}} \mathrm{d}a_{\circledcirc} +
		\sum_{\circledcirc \in \{\MA, \VE\}} \int_{A_{\circledcirc}} m_{\circledast \circledcirc}(a_{\circledcirc}) \int_{0}^{a_{\circledcirc}} e^{-\beta(a_{\circledcirc}
		- \alpha_{\circledcirc})} v_{\circledcirc}(\alpha_{\circledcirc}) \mathrm{d}\alpha_{\circledcirc} \mathrm{d}a_{\circledcirc}.
		\notag
	\end{equation}
	The latter can be written as
	\begin{equation}
		\sum_{\circledcirc \in \{\MA, \VE\}} M_{\circledast \circledcirc}(\beta) c_{\circledcirc} =
		b_{\circledcirc}(\beta) \text{ for } \circledast \in \{\MA, \VE\} \notag
	\end{equation}
	with
	\begin{equation}
		\begin{split}
			M_{\circledast \circledcirc}(\beta) &:= \delta_{\circledast \circledcirc} -
			\sum_{\circledcirc \in \{\MA, \VE\}} \int_{A_{\circledcirc}} m_{\circledast \circledcirc}(a_{\circledcirc}) e^{-\beta a_{\circledcirc}} \mathrm{d}a_{\circledcirc} \text{ for } \circledast, \circledcirc \in \{\MA, \VE\}, \\
			b_{\circledast}(\beta) &:= \sum_{\circledcirc \in \{\MA, \VE\}} \int_{A_{\circledcirc}} m_{\circledast \circledcirc}(a_{\circledcirc}) \int_{0}^{a_{\circledast}} e^{-\beta(a_{\circledast} - \alpha_{\circledast})}
			v_{\circledast}(\alpha_{\circledast}) \mathrm{d}\alpha_{\circledast} \mathrm{d}a_{\circledcirc} \text{ for } \circledast \in \{\MA, \VE\}.
		\end{split}
		\notag
	\end{equation}
	Since we can estimate
	\begin{equation}
		\begin{split}
			|M_{\circledast \circledcirc}(\beta) - \delta_{\circledast \circledcirc}| &\leq
			\sum_{\circledcirc \in \{\MA, \VE\}} \int_{A_{\circledcirc}} |m_{\circledast \circledcirc}(a_{\circledcirc}) e^{-\beta a_{\circledcirc}}| \mathrm{d}a_{\circledcirc} \\
			&\leq \sum_{\circledcirc \in \{\MA, \VE\}} \|m_{\circledast \circledcirc}\|_{L^{\infty}(A_{\circledcirc})} \tfrac{1 - e^{-\beta a_{\circledcirc}^{\dag}}}{\beta}
			(a_{\circledcirc}^{\dag})^{1/2} \to 0 \text{ for } \beta \to \infty,
		\end{split}
		\notag
	\end{equation}
	there exists a number $\beta_{1} > 0$ such that the matrix $M(\beta) \in \RR^{2 \times 2}$ is invertible for all $\beta > \beta_{1}$
	with its inverse matrix given as a Neumann series.
	Further, the vector $b(\beta) \in \RR^{2}$ is well-defined since
	\begin{equation}
		\begin{split}
			|b_{\circledast}| &\leq
			\sum_{\circledcirc \in \{\MA, \VE\}} \int_{A_{\circledcirc}} |m_{\circledast \circledcirc}(a_{\circledcirc})| \int_{0}^{a_{\circledast}}
			|v_{\circledast}(\alpha_{\circledast})| \mathrm{d}\alpha_{\circledast} \mathrm{d}a_{\circledcirc} \\
			&\leq \sum_{\circledcirc \in \{\MA, \VE\}} \|m_{\circledast \circledcirc}\|_{L^{\infty}(A_{\circledcirc})} (a_{\circledcirc}^{\dag})^{3/2} \|v_{\circledcirc}\|_{L^{2}(A_{\circledcirc})} < \infty.
		\end{split}
		\notag
	\end{equation}
	Moreover, we see that the expression $\|b(\beta)\|_{\RR^{2}}$ linearly depends on $\|v\|_{X}$
	whereas $\|M^{-1}(\beta)\|_{\RR^{2 \times 2}}$ does not depend on $v$.
	Therefore,
	\begin{equation}
		\|c\|_{\RR^{2}} \leq C_{2}(\beta) \|v\|_{X} \text{ for some } C_{2}(\beta) > 0. \notag
	\end{equation}
	Plugging this into Equation (\ref{EQUATION_CONSTRUCTION_OF_U_FOR_BABUSKA_LAX_MILGRAM_SOLUTION}),
	we obtain a solution $u \in H^{1}(A_{\MA}) \times H^{1}(A_{\VE})$ satisfying Equations (\ref{EQUATION_CONSTRUCTION_OF_U_FOR_BABUSKA_LAX_MILGRAM}),
	(\ref{EQUATION_CONSTRUCTION_OF_U_FOR_BABUSKA_LAX_MILGRAM_IC}) and thus lying in $D(A)$.
	By construction, we obtain
	\begin{equation}
		a(u, v; \beta) = \sum_{\circledast \in \{\MA, \VE\}} \langle \partial_{a_{\circledast}} u_{\circledast} + \beta u_{\circledast}, v_{\circledast}\rangle_{L^{2}(A_{\circledast})}
		= \sum_{\circledast \in \{\MA, \VE\}} \|v_{\circledast}\|_{L^{2}(A_{\circledast})}^{2} = \|v\|_{X}^{2} \notag
	\end{equation}
	and
	\begin{equation}
		\|u\|_{D(A)} = \|u\|_{H^{1}(A_{\circledast}) \times H^{1}(A_{\circledcirc})} \leq C_{1}(\beta)(1 + C_{2}(\beta)) \|v\|_{X}. \notag
	\end{equation}
	Thus, the bilinear form $a$ satisfies the inf-sup-condition
	meaning that the operator $A - \beta \mathrm{id}_{X}$ is continuously invertible and therefore surjective.
	
	Altogether we have shown that $A - \beta \mathrm{id}_{X}$ is m-dissipative for $\beta \geq \max\{\beta_{0}, \beta_{1}\}$.
\end{proof}

Taking into account Lemmas \ref{LEMMA_OPERATOR_A_DENSE_AND_CLOSED} and \ref{LEMMA_A_MINUS_BETA_M_DISSIPATIVE},
we apply the Theorem of Lumer \& Phillips
as well as the well-known perturbation result for bounded operators (cf. \cite[Corollary 1.3]{Pa1983}) to conclude
\begin{theorem}
	The operator $A$ is a generator of a $C_{0}$-semigroup $(S(t))_{t \geq 0}$ of bounded linear operators on $X$.
\end{theorem}

Now, we exploit \cite[Theorem 3.1.12]{ArBaHieNeu2001} and \cite[Corollary 3.1.17]{ArBaHieNeu2001} and conclude
\begin{theorem}
	\label{THEOREM_STRONG_SOLUTION_EXISTENCE_AND_UNIQUENESS}
	Assume that $u^{0} := (u_{\MA}^{0}, u_{\VE}^{0}) \in X$, $f \in L^{2}(0, T; X)$.
	Then there exists a unique mild solution $u \in C^{0}([0, \infty), X)$ to Equation (\ref{EQUATION_OPERATOR_EQUATION}) given as
	\begin{equation}
		u(t) = S(t) u^{0} + \int_{0}^{t} S(t - s) f(s) \mathrm{d}s \text{ for } t \geq 0 \notag
	\end{equation}
	continuously depending on the data in sense of the existence of constants $C \geq 1$, $\omega \in \RR$ such that
	\begin{equation}
		\|u\|_{L^{\infty}(0, T; X)} 
		\leq C \big(1 + e^{\omega T}\big) \big(\|u^{0}\|_{X} + \|f\|_{L^{2}(0, T; X)}\big) \text{ for any } T > 0. \notag
	\end{equation}
	If $u^{0} \in D(A)$ and $f \in H^{1}(0, T; X)$,
	then there exists a constant $C > 0$ such that Equation (\ref{EQUATION_MODEL_ABSTRACT_FORM})
	possesses a unique classical solution
	\begin{equation}
		u \in C^{1}([0, T], X) \cap C^{0}([0, T], D(A)). \notag
	\end{equation}
\end{theorem}

Finally, we want to study the asymptotic behavior of solutions to (\ref{EQUATION_MODEL_EQUATION_1})--(\ref{EQUATION_MODEL_EQUATION_6})
in the absense of immigration or emigration, i.e., $f \equiv 0_{X}$.
We define the ``natural'' energy via
\begin{equation}
	E(t) := \tfrac{1}{2} \Big(\int_{A_{\MA}} u_{\MA}^{2}(t, a_{\MA}) \mathrm{d}a_{\MA} + \int_{A_{\VE}} u_{\VE}^{2}(t, a_{\VE}) \mathrm{d}a_{\VE}\Big)
	= \tfrac{1}{2} \|u(t)\|_{X}^{2} \notag
\end{equation}
and easily see that the exponential stability of the zero solution to (\ref{EQUATION_MODEL_EQUATION_1})--(\ref{EQUATION_MODEL_EQUATION_6})
is equivalent with the exponential stability of the zero solution to (\ref{EQUATION_MODEL_ABSTRACT_FORM})
whereas the latter holds true if and only if the semigroup $(S(t))_{t \geq 0}$ is exponentially stable.

\begin{theorem}
	Assume that
	\begin{equation}
		\max_{\circledast \in \{\MA, \VE\}} \sum_{\circledcirc \in \{\MA, \VE\}} a_{\circledast}^{\dag} a_{\circledcirc}^{\dag} \|m_{\circledast \circledcirc}\|_{L^{\infty}(A_{\circledcirc})}^{2} < \tfrac{1}{4}.
		\notag
	\end{equation}
	Then the energy $E(t)$ decays exponentially to zero for $t \to \infty$, i.e.,
	\begin{equation}
		E(t) \leq C e^{-2 \alpha t} E(0) \text{ for } t \geq 0 \notag
	\end{equation}
	with
	\begin{equation}
		\alpha :=
		\frac{\min
		\Big\{1 - 4 \sum\limits^{}_{\circledcirc \in \{\MA, \VE\}} a_{\circledast}^{\dag} a_{\circledcirc}^{\dag} \|m_{\circledast \circledcirc}\|_{L^{\infty}(A_{\circledcirc})}^{2} \,\big|\, \circledast \in \{\MA, \VE\}\Big\}}
		{2 \max\{a_{\MA}^{\dag}, a_{\VE}^{\dag}\}}, \quad
		C := \frac{1}{2 \max\{a_{\MA}^{\dag}, a_{\VE}^{\dag}\}}. \notag
	\end{equation}
\end{theorem}

\begin{proof}
	Since any initial data $u^{0} \in X$ can be approximated by a sequence from $D(A)$,
	we assume without loss of generality that $u^{0} \in D(A)$
	and denote by $u$ the corresponding unique classical solution of Equation (\ref{EQUATION_MODEL_ABSTRACT_FORM}),
	which in its turn is a classical solution to (\ref{EQUATION_MODEL_EQUATION_1})--(\ref{EQUATION_MODEL_EQUATION_6}).

	We consider the Lyapunov functional
	\begin{equation}
		F(t; u) := \sum_{\circledast \in \{\MA, \VE\}} \int_{A_{\circledast}} (2 a_{\circledast}^{\dag} - a_{\circledast}) u_{\circledast}^{2}(t, a_{\circledast}) \mathrm{d}a_{\circledast}. \notag
	\end{equation}
	Obviously,
	\begin{equation}
		0 \leq E(t; u) \leq F(t; u) \leq 2 \max\{a_{\MA}^{\dag}, a_{\VE}^{\dag}\} E(t; u). \notag
	\end{equation}
	Moreover, $F(\cdot; u)$ is Frech\'{e}t differentiable along the solution $u$
	and due to Equations (\ref{EQUATION_MODEL_TRANSFORMED_EQUATION_1})--(\ref{EQUATION_MODEL_TRANSFORMED_EQUATION_4}) satisfies
	\begin{align*}
		\frac{\mathrm{d}}{\mathrm{d}t} F(t; u) &=
		\sum_{\circledast \in \{\MA, \VE\}} \int_{A_{\circledast}} (2 a_{\circledast}^{\dag} - a_{\circledast}) \partial_{t} u_{\circledast}(t, a_{\circledast}) u_{\circledast}(t, a_{\circledast}) \mathrm{d}a_{\circledast} \\
		&=
		-\sum_{\circledast \in \{\MA, \VE\}} \int_{A_{\circledast}} (2 a_{\circledast}^{\dag} - a_{\circledast}) \partial_{a_{\circledast}} u_{\circledast}(t, a_{\circledast}) u_{\circledast}(t, a_{\circledast}) \mathrm{d}a_{\circledast} \\
		&=
		-\tfrac{1}{2} \sum_{\circledast \in \{\MA, \VE\}} \int_{A_{\circledast}} (2 a_{\circledast}^{\dag} - a_{\circledast}) \partial_{a_{\circledast}} \big(u_{\circledast}^{2}(t, a_{\circledast})\big) \mathrm{d}a_{\circledast} \\
		&=
		-\tfrac{1}{2} \sum_{\circledast \in \{\MA, \VE\}}
		\Big(\int_{A_{\circledast}} u_{\circledast}^{2}(t, a_{\circledast}) \mathrm{d}a_{\circledast} +
		(2 a_{\circledast}^{\dag} - a_{\circledast}) u_{\circledast}^{2}(t, a_{\circledast})\Big|_{a_{\circledast} = 0}^{a_{\circledast} = a_{\circledast}^{\dag}}\Big) \\ 
		&\leq
		-\tfrac{1}{2} \sum_{\circledast \in \{\MA, \VE\}}
		\Big(\int_{A_{\circledast}} u_{\circledast}^{2}(t, a_{\circledast}) \mathrm{d}a_{\circledast} - 2 a_{\circledast} u_{\circledast}^{2}(t, 0)\Big) \\
		&\leq
		-\tfrac{1}{2} \sum_{\circledast \in \{\MA, \VE\}}
		\Big[\int_{A_{\circledast}} u_{\circledast}^{2}(t, a_{\circledast}) \mathrm{d}a_{\circledast} -
		2 a_{\circledast} \Big(\sum_{\circledcirc \in \{\MA, \VE\}} \int_{A_{\circledcirc}} m_{\circledast \circledcirc}(a_{\circledcirc}) u_{\circledcirc}(t, a_{\circledcirc}) \mathrm{d}a_{\circledcirc}\Big)^{2}\Big] \\
		&\leq
		-\tfrac{1}{2} \sum_{\circledast \in \{\MA, \VE\}}
		\Big[\int_{A_{\circledast}} u_{\circledast}^{2}(t, a_{\circledast}) \mathrm{d}a_{\circledast} \\
		&-4 a_{\circledast} \sum_{\circledcirc \in \{\MA, \VE\}} a_{\circledcirc} \|m_{\circledast \circledcirc}\|_{L^{\infty}(A_{\circledcirc})}^{2} \int_{A_{\circledcirc}} u_{\circledcirc}^{2}(t, a_{\circledcirc}) \mathrm{d}a_{\circledcirc}\Big] \\
		&= -\sum_{\circledast \in \{\MA, \VE\}}
		\Big(1 - 4 \sum_{\circledcirc \in \{\MA, \VE\}} a_{\circledast}^{\dag} a_{\circledcirc}^{\dag} \|m_{\circledast \circledcirc}\|_{L^{\infty}(A_{\circledcirc})}^{2}\Big) \int_{A_{\circledast}} u_{\circledast}^{2}(t, a_{\circledast}) \mathrm{d}a_{\circledast} \\
		&= -\min\Big\{1 - 4 \sum_{\circledcirc \in \{\MA, \VE\}} a_{\circledast}^{\dag} a_{\circledcirc}^{\dag} \|m_{\circledast \circledcirc}\|_{L^{\infty}(A_{\circledcirc})}^{2} \,\big|\, \circledast \in \{\MA, \VE\}\Big\} E(t, u)
		= -2\alpha F(t, u),
	\end{align*}
	where we performed an integration by parts and used Young's and H\"older's inequalities.
	Now, applying Gronwall's inequality, we obtain
	\begin{equation}
		E(t, u) \leq F(t, u) \leq e^{-2\alpha t} F(0, u) \leq C e^{-2\alpha t} E(0, u) \text{ for } t \geq 0, \notag
	\end{equation}
	which was our claim.
\end{proof}

\section{Finite Difference Scheme and Convergence Analysis}
\label{SECTION_FINITE_DIFFERENCE_SCHEME}
In this section, we propose an implicit finite difference method
to numerically solve the initial-boundary value problem (\ref{EQUATION_MODEL_EQUATION_1})--(\ref{EQUATION_MODEL_EQUATION_6}).
Under minimal regularity assumptions on the data, we show the scheme to be convergent.
In our investigations, we decided to depart from the standard approach of assuming the $C^{2}$-differentiability of solutions (cf., e.g., \cite{ArbMil1989}),
since, to assure for this high regularity of solutions,
one would require in addition to an extra smoothness condition on the data and system parameters some rather restrictive compatibility conditions on $u^{0}$ and $f$
which are usually not satisfied in real applications.
Though finite difference discretizations of Equations (\ref{EQUATION_MODEL_EQUATION_1})--(\ref{EQUATION_MODEL_EQUATION_6})
satisfy the Courant-Friedrichs-Levy condition,
we decided to use an implicit scheme instead of an explicit one
to assure for better stability on long time horizons.
To the authors' best knowledge, earlier works (viz. \cite{AbLo1997}, \cite{Su1994}, etc.) do not provide
a rigorous convergence study for the implicit scheme in $L^{2}$-settings, in particular, under minimal regularity assumptions.
For studies on explicit schemes we refer the reader to \cite{ArbMil1989}, \cite{Ko2003}.

Throughout this section, we assume that $m_{\circledast \circledcirc} \in C^{0}(\bar{A}_{\circledcirc})$ for $\circledast, \circledcirc \in \{\MA, \VE\}$ and
\begin{equation}
	u^{0} \in D(A), \quad f \in H^{1}(0, T; X) \cap C^{0}\big([0, T], C^{0}(\bar{A}_{\MA}) \times C^{0}(\bar{A}_{\VE})\big). \notag
\end{equation}
Then, the conditions of Theorem \ref{THEOREM_STRONG_SOLUTION_EXISTENCE_AND_UNIQUENESS} are trivially fulfilled and
we obtain a unique strong solution of Equation (\ref{EQUATION_MODEL_ABSTRACT_FORM}).
Again, it should be stressed that no compatibility conditions are required here.

Selecting the age discretization steps
\begin{equation}
	h_{\circledast} = a_{\circledast}^{\dag}/N_{\circledast, h_{\circledast}} \text{ such that }
	N_{\circledast, h_{\circledast}} \in \mathbb{N} \text{ for } \circledast \in \{\MA, \VE\}, \notag
\end{equation}
we define the equidistant age lattices
\begin{equation}
	A_{\circledast}^{h} := \{a^{h_{\circledast}}_{\circledast, i_{\circledast}} \,|\, i_{\circledast} = 0, \dots, N_{\circledast, h_{\circledast}}\}
	\text{ with } a^{h_{\circledast}}_{\circledast, i_{\circledast}} := i_{\circledast} h_{\circledast} \text{ for } i_{\circledast} = 0, \dots, N_{\circledast, h_{\circledast}}
	 \notag
\end{equation}
as well as their ``interiors'' and ``boundaries''
\begin{equation}
	\stackrel{\circ}{A_{\circledast}^{h_{\circledast}}} := \{a^{h_{\circledast}}_{\circledast, i_{\circledast}} \,|\, i_{\circledast} = 1, \dots, N_{\circledast, h_{\circledast}}\}
	\text{ and }
	\partial A_{\circledast}^{h_{\circledast}} := \{a^{h_{\circledast}}_{\circledast, 0}\} = \{0\}, \text{ respectively}, \notag
\end{equation}
for $\circledast \in \{\MA, \VE\}$.
In this section, we adopt the notation from the Appendix letting
$L_{\tau}^{2}, L_{h_{\MA}}^{2}, L_{h_{\VE}}^{2}$ denote discrete Lebesgue spaces.

For each time $t \geq 0$, the functions $u_{\circledast}(t, \cdot)$ and $f_{\circledast}(t, \cdot)$ will be approximated by the lattice functions
$u_{\circledast}^{h_{\circledast}}(t, \cdot), f_{\circledast}^{h_{\circledast}}(t, \cdot) \colon A_{\circledast}^{h_{\circledast}} \to \RR$ for $\circledast \in \{\MA, \VE\}$.
Using the backwards difference approximation for the age derivatives and a Riemann sum discretization for the integral,
we obtain the following semi-discretization with respect to the age variables
\begin{align}
	\partial_{t} u_{\circledast, i_{\circledast}}^{h_{\circledast}}(t) + \tfrac{u_{\circledast, i_{\circledast}}^{h_{\circledast}}(t) - u_{\circledast, i_{\circledast} - 1}^{h_{\circledast}}(t)}{h_{\circledast}} =
	f_{\circledast, i_{\circledast}}^{h_{\circledast}}(t)
	\text{ for } i_{\circledast} = 1, \dots, N_{\circledast, h_{\circledast}}, \circledast \in \{\MA, \VE\}, t > 0,
	\label{EQUATION_MODEL_TRANSFORMED_DISCRETIZED_EQUATION_1} \\
	u_{\circledast, 0}^{h_{\circledast}}(t) = \sum\limits^{}_{\circledcirc \in \{\MA, \VE\}} h_{\circledcirc} \sum_{i_{\circledcirc} = 1}^{N_{\circledcirc, h_{\circledast}}} m_{\circledast \circledcirc}(a_{\circledcirc, i_{\circledcirc}}^{h_{\circledcirc}}) u_{\circledcirc, i_{\circledcirc}}^{h_{\circledcirc}}(t)
	\text{ for } \circledast \in \{\MA, \VE\}, t > 0,
	\label{EQUATION_MODEL_TRANSFORMED_DISCRETIZED_EQUATION_2} \\
	u_{\circledast, i_{\circledast}}^{h_{\circledast}}(0) =
	u_{\circledast, i_{\circledast}}^{0, h_{\circledast}}
	\text{ for } i_{\circledast} = 1, \dots, N_{\circledast, h_{\circledast}}, \circledast \in \{\MA, \VE\}
	\label{EQUATION_MODEL_TRANSFORMED_DISCRETIZED_EQUATION_3}
\end{align}
with $u_{\circledast, i_{\circledast}}^{h_{\circledast}}(t)$ and $f_{\circledast, i_{\circledast}}^{h_{\circledast}}$
approximating $u_{\circledast}(t, a^{h_{\circledast}}_{\circledast, i_{\circledast}})$ and $f_{\circledast}(t, a^{h_{\circledast}}_{\circledast, i_{\circledast}})$, respectively,
and $u_{\circledast, i_{\circledast}}^{0, h_{\circledast}}$ being an approximation for $u_{\circledast}^{0}(a_{\circledast, i_{\circledast}}^{h_{\circledast}})$ for $\circledast \in \{\MA, \VE\}$.

We let
\begin{equation}
	X^{h} := L^{2}_{h_{\MA}}(A_{\MA}^{h_{\MA}}) \times L^{2}_{h_{\VE}}(A_{\VE}^{h_{\VE}}), \quad
	\stackrel{\circ}{X^{h}} := L_{h_{\MA}}^{2}(\stackrel{\circ}{A_{\MA}^{h_{\MA}}}) \times L_{h_{\VE}}^{2}(\stackrel{\circ}{A_{\VE}^{h_{\VE}}}) \notag
\end{equation}
and define the restriction operators
\begin{equation}
	\stackrel{\circ}{u^{h}_{\circledast}} := u^{h}_{\circledast}(t, \cdot)\big|_{\stackrel{\circ}{A_{\circledast}^{h_{\circledast}}}}, \quad
	\stackrel{\diamond}{u^{h}_{\circledast}} := u^{h}_{\circledast}(t, \cdot)\big|_{{\partial A_{\circledast}^{h_{\circledast}}}}
	\text{ for } u^{h} \in X^{h}, \circledast \in \{\MA, \VE\}. \notag
\end{equation}
Further, we introduce the linear operators
$\stackrel{\circ}{B^{h}} \colon \stackrel{\circ}{X^{h}} \to \RR^{2}$ and $\stackrel{\circ}{A^{h}} \colon D(\stackrel{\circ}{A^{h}}) \to \stackrel{\circ}{X^{h}}$ by the means of
\begin{align*}
	\stackrel{\circ}{B^{h}_{\circledast}} \stackrel{\circ}{u^{h}} &=
	\sum\limits_{\circledcirc \in \{\MA, \VE\}} h_{\circledcirc} \sum\limits_{i_{\circledcirc} = 1}^{N_{\circledcirc, h_{\circledcirc}}} m_{\circledast \circledcirc}(a_{\circledcirc, i_{\circledcirc}}^{h_{\circledcirc}}) \stackrel{\circ}{u}_{\circledcirc, i_{\circledcirc}}^{h_{\circledcirc}}
	\text{ for } \circledast \in \{\MA, \VE\}, \\
	\Big[\stackrel{\circ}{A^{h}_{\circledast}} \stackrel{\circ}{u^{h}}\Big](a^{h_{\circledast}}_{h_{\circledast}, i_{\circledast}}) &=
	\begin{cases}
		-\tfrac{\stackrel{\circ}{u^{h}}_{\circledast, 1}^{h_{\circledast}} - \stackrel{\circ}{B^{h}_{\circledast}} \stackrel{\circ}{u^{h}}}{h_{\circledast}}, & i_{\circledast} = 1 \\
		-\tfrac{\stackrel{\circ}{u^{h}}_{\circledast, i_{\circledast}}^{h_{\circledast}} - \stackrel{\circ}{u}_{\circledast, i_{\circledast} - 1}^{h_{\circledast}}}{h_{\circledast}}, & 
		i_{\circledast} \in \{2, \dots, N_{\circledast, h_{\circledast}}\}
	\end{cases}
	\text{ for } \circledast \in \{\MA, \VE\}
	\notag
\end{align*}
where $D(\stackrel{\circ}{A^{h}}) := \stackrel{\circ}{X^{h}}$ is equipped with the inner product
\begin{equation}
	\langle u^{h}, v^{h}\rangle_{D(\stackrel{\circ}{A^{h}})} :=
	\langle u^{h}, v^{h}\rangle_{\stackrel{\circ}{X^{h}}} +
	\langle \stackrel{\circ}{A^{h}} u^{h}, \stackrel{\circ}{A^{h}} v^{h}\rangle_{\stackrel{\circ}{X^{h}}} \text{ for } u^{h}, v^{h} \in D(\stackrel{\circ}{A^{h}}).
	\notag
\end{equation}
Hence, Equations (\ref{EQUATION_MODEL_TRANSFORMED_DISCRETIZED_EQUATION_1})--(\ref{EQUATION_MODEL_TRANSFORMED_DISCRETIZED_EQUATION_3}) can be equivalently transformed to
\begin{align}
	\partial_{t} \stackrel{\circ}{u^{h}}(t) &= \; \stackrel{\circ}{A^{h}} \stackrel{\circ}{u^{h}}(t) + \stackrel{\circ}{f^{h}}(t) \text{ for } t > 0, \label{EQUATION_MODEL_TRANSFORMED_DISCRETIZED_EQUATION_ODE_1} \\
	\stackrel{\diamond}{u^{h}}(t) &= B^{h} \stackrel{\circ}{u^{h}}(t) \text{ for } t > 0, \label{EQUATION_MODEL_TRANSFORMED_DISCRETIZED_EQUATION_ODE_2} \\
	u^{h}(0) &= u^{0, h} \label{EQUATION_MODEL_TRANSFORMED_DISCRETIZED_EQUATION_ODE_3}
\end{align}
where $u^{0, h}$ and $f^{h, \tau}$ are approximations of $u^{0}$ and $f$, respectively.

For $T > 0$, we consider a time step $\tau = \tfrac{T}{N_{\tau}}$ with $N_{\tau} \in \NN$ and define the time lattice
\begin{equation}
	Z^{\tau} := \{t_{k}^{\tau} \,|\, k = 0, \dots, N^{\tau}\} \text{ with } t_{k}^{\tau} := \tau k \text{ for } k = 0, \dots, N^{\tau} \notag
\end{equation}
as well as its ``interior'' $\stackrel{\circ}{Z^{\tau}} := \{t_{k}^{\tau} \,|\, k = 1, \dots, M\}$.
The functions $u^{h}, f^{h} \colon [0, T] \to X^{h}$ will now be approximated by the lattice functions
$u^{h, \tau}, f^{h, \tau} \colon Z_{\tau} \to X^{h}$.
Similarly, $\stackrel{\circ}{u^{h}}$, $\stackrel{\circ}{f^{h}} \colon [0, T] \to \stackrel{\circ}{X^{h}}$
will be approximated by $\stackrel{\circ}{u^{h, \tau}}, \stackrel{\circ}{f^{h, \tau}} \colon Z^{\tau} \to \stackrel{\circ}{X^{h}}$.

For $\vartheta \in [0, 1]$, the ODE system (\ref{EQUATION_MODEL_TRANSFORMED_DISCRETIZED_EQUATION_ODE_1}), (\ref{EQUATION_MODEL_TRANSFORMED_DISCRETIZED_EQUATION_ODE_3}) can now be discretized using the $\vartheta$-method
whereas Equation (\ref{EQUATION_MODEL_TRANSFORMED_DISCRETIZED_EQUATION_ODE_2}) will just be restricted onto the inner time grid $\stackrel{\circ}{Z^{\tau}}$.
This yields a difference equation for $u^{h, \tau}$
\begin{align}
	\tfrac{\stackrel{\circ}{u^{h, \tau}}(t_{k}, \cdot) - \stackrel{\circ}{u^{h, \tau}}(t_{k-1}, \cdot)}{\tau} &= \vartheta \stackrel{\circ}{A^{h}} \stackrel{\circ}{u^{h, \tau}}(t_{k}, \cdot) + (1 - \vartheta) \stackrel{\circ}{A^{h}} \stackrel{\circ}{u^{h, \tau}}(t_{k - 1}, \cdot) \label{EQUATION_FULL_DISCRETIZATION_1} \\
	&+ \vartheta \stackrel{\circ}{f^{h, \tau}}(t_{k}, \cdot) + (1 - \vartheta) \stackrel{\circ}{f^{h, \tau}}(t_{k-1}, \cdot) \text{ for } k = 1, \dots, N^{\tau}, \notag \\
	\stackrel{\diamond}{u^{h, \tau}}(t_{k}, \cdot) &= \; \stackrel{\circ}{B^{h}} \stackrel{\circ}{u^{h, \tau}}(t_{k}, \cdot) \text{ for } k = 1, \dots, N^{\tau}, \label{EQUATION_FULL_DISCRETIZATION_2} \\
	u^{h, \tau}(0, \cdot) &= u^{0, h}. \label{EQUATION_FULL_DISCRETIZATION_3}
\end{align}

Next, we define the bounded linear operators
\begin{align*}
	\mathcal{L}^{h, \tau} \colon D(\mathcal{L}^{h, \tau}) &\to
	X^{h} \times L^{2}_{\tau}(\stackrel{\circ}{{Z}^{\tau}}, \stackrel{\circ}{X^{h}}) \times L^{2}_{\tau}(\stackrel{\circ}{{Z}^{\tau}}, \RR^{2}), \\
	u^{h, \tau}
	&\mapsto {\scriptsize
	\begin{pmatrix}
		u^{h, \tau}(0, \cdot) \\
		\tfrac{u^{h, \tau}(t_{k}^{\tau}, \cdot) - u^{h, \tau}(t_{k-1}^{\tau}, \cdot)}{\tau} - \vartheta (\stackrel{\circ}{A^{h}} u^{h, \tau})(t_{k}^{\tau}, \cdot) -
		(1 - \vartheta) (\stackrel{\circ}{A^{h}} u^{h, \tau})(t_{k-1}^{\tau}, \cdot), \quad
		k = 1, \dots, N^{\tau} \\
		\stackrel{\diamond}{u^{h, \tau}}(t_{k}^{\tau}, \cdot) - \stackrel{\circ}{B^{h}} u^{h, \tau}(t_{k}^{\tau}, \cdot), \quad k = 1, \dots, N^{\tau}
	\end{pmatrix}}
	\notag
\end{align*}
with
\begin{equation}
	D(\mathcal{L}^{h, \tau}) :=
	H^{1}_{\tau}\big(Z^{\tau}, H^{1}_{h_{\MA}}(A_{\MA}^{h_{\MA}}) \times H^{1}_{h_{\VE}}(A_{\VE}^{h_{\VE}})\big) \notag
\end{equation}
and
\begin{align*}
	\mathcal{F}^{h, \tau} \colon X^{h} \times L^{2}_{\tau}(Z^{\tau}, \stackrel{\circ}{X^{h}}) &\to
	X^{h} \times L^{2}_{\tau}(\stackrel{\circ}{Z^{\tau}}, \stackrel{\circ}{X^{h}}) \times L^{2}_{\tau}(\stackrel{\circ}{{Z}^{\tau}}, \RR^{2}), \\
	(u^{0, h}, \stackrel{\circ}{f^{h, \tau}}) &\mapsto {\scriptsize
	\begin{pmatrix}
		u^{0, h} \\
		\vartheta \stackrel{\circ}{f^{h, \tau}}(t_{k}^{\tau}, \cdot) + (1 - \vartheta) \stackrel{\circ}{f^{h, \tau}}(t_{k-1}^{\tau}, \cdot), \quad k = 1, \dots, N^{\tau} \\
		0, \quad k = 1, \dots, N^{\tau}
	\end{pmatrix}.}
\end{align*}
With this notation, Equations (\ref{EQUATION_FULL_DISCRETIZATION_1})--(\ref{EQUATION_FULL_DISCRETIZATION_3}) can be equivalently re-written as
\begin{equation}
	\mathcal{L}^{h, \tau} u^{\tau, h} = \mathcal{F}^{h, \tau} (u^{0, h}, \stackrel{\circ}{f^{h, \tau}}).
	\label{EQUATION_NUMERICAL_SCHEME_ABSTRACT_FORM}
\end{equation}
Investigating the solvability of the numerical scheme (\ref{EQUATION_NUMERICAL_SCHEME_ABSTRACT_FORM})
as well as its convergence for $(h, \tau) \to 0$ will be our thrust for the rest of this section.

\subsection{Consistency}
To prove the consistency for the difference scheme (\ref{EQUATION_NUMERICAL_SCHEME_ABSTRACT_FORM}),
we exploit basic approximation properties of Banach space-valued functions and Bochner integrals (see, e.g., \cite[Chapter 1]{ArBaHieNeu2001}).
No error estimates based on Taylor expansion will be used here due to the possible lack of classical differentiability
in real-world applications.

By the virtue of Sobolev embedding theorem, we have
\begin{equation}
	D(A) \hookrightarrow H^{1}(A_{\MA}) \times H^{1}(A_{\VE}) \hookrightarrow C^{0}(\bar{A}_{\MA}) \times C^{0}(\bar{A}_{\VE}) \notag
\end{equation}
as well as
\begin{equation}
	C^{0}([0, T], D(A)) \hookrightarrow C^{0}\big([0, T], C^{0}(\bar{A}_{\MA}) \times C^{0}(\bar{A}_{\VE})\big). \notag
\end{equation}
Hence, the elements from $D(A)$ and $C^{0}([0, T], D(A))$, being in general some Lebesgue equivalence classes,
have a continuous representative and thus can be evaluated pointwise.

\begin{lemma}
	\label{LEMMA_CONSISTENCY_FOR_A_H}
	For $u \in X$ let $u^{h} := \Big(u_{\circledast}\big|_{A_{\circledast}^{h_{\circledast}}}\Big)_{\circledast \in \{\MA, \VE\}}$.
	Then \\[-0.5cm]
	\begin{enumerate}
		\item[{\it i)}]
		$\sum\limits^{}_{\circledast \in \{\MA, \VE\}} \Big|\stackrel{\circ}{B_{\circledast}^{h}} \stackrel{\circ}{u^{h}} -
		\sum\limits^{}_{\circledcirc \in \{\MA, \VE\}} \int_{A_{\circledcirc}} m_{\circledast \circledcirc}(a_{\circledcirc}) m_{\circledcirc}(a_{\circledcirc}) \mathrm{d}a_{\circledcirc}\Big|^{2} \to 0$ as $h \to 0$.
	
		\item[{\it ii)}]
		$\sum\limits^{}_{\circledast \in \{\MA, \VE\}} \sum\limits_{i_{\circledast} = 1}^{N_{\circledast, h_{\circledast}}}
		\int\limits_{a_{\circledast, i_{\circledast} - 1}^{h_{\circledast}}}^{a_{\circledast, i_{\circledast}}^{h_{\circledast}}}
		\Big(\big[\stackrel{\circ}{A^{h}} u^{h}\big]_{\circledast}(a_{\circledast, i_{\circledast}}^{h_{\circledast}}) - A_{\circledast}(a_{\circledast})\Big)^{2} \mathrm{d}a_{\circledast} \to 0$ as $h \to 0$ if $u \in D(A)$.
	\end{enumerate}
\end{lemma}

\begin{proof}
	\begin{enumerate}
		\item[{\it i)}] Using Lemma \ref{LEMMA_APPENDIX_INTEGRATION_OPERATOR_ESTIMATE}, we trivially obtain
		\begin{align*}
			&\sum_{\circledast \in \{\MA, \VE\}} \Big|\stackrel{\circ}{B_{\circledast}^{h}} \stackrel{\circ}{u^{h}} -
			\sum_{\circledcirc \in \{\MA, \VE\}} \int_{A_{\circledcirc}} m_{\circledast \circledcirc}(a_{\circledcirc}) m_{\circledcirc}(a_{\circledcirc}) \mathrm{d}a_{\circledcirc}\Big|^{2} \\
			&\leq
			\sum_{\circledcirc, \circledast \in \{\MA, \VE\}} h_{\circledcirc}^{-1} \Big|\int_{A_{\circledcirc}} m_{\circledast\circledcirc}(a_{\circledcirc}) u_{\circledcirc}(a_{\circledcirc}) \mathrm{d}a_{\circledcirc} -
			\sum_{i_{\circledcirc} = 1}^{N_{\circledcirc, h_{\circledcirc}}} m_{\circledast \circledcirc}(a_{\circledcirc, i_{\circledcirc}}^{h_{\circledcirc}}) u_{\circledcirc}(a_{\circledcirc, i_{\circledcirc}}^{h_{\circledcirc}}) \Big|^{2}
			\to 0
		\end{align*}
		as $h \to 0$.

		\item[{\it ii)}] Using {\it i)} as well as Lemma \ref{LEMMA_APPENDIX_DIFFERENCE_OPERATOR_ESTIMATE} and applying Young's inequality, we get
		\begin{align*}
			&\sum\limits^{}_{\circledast \in \{\MA, \VE\}} \sum_{i_{\circledast} = 1}^{N_{\circledast, h_{\circledast}}}
			\int_{a_{\circledast, i_{\circledast} - 1}^{h_{\circledast}}}^{a_{\circledast, i_{\circledast}}^{h_{\circledast}}}
			\Big(\big[\stackrel{\circ}{A^{h}} u^{h}\big]_{\circledast}(a_{\circledast, i_{\circledast}}^{h_{\circledast}}) - A_{\circledast}(a_{\circledast})\Big)^{2} \mathrm{d}a_{\circledast} \\
			&= \sum\limits^{}_{\circledast \in \{\MA, \VE\}}
			\int_{a_{\circledast, 0}^{h_{\circledast}}}^{a_{\circledast, 1}^{h_{\circledast}}} \Big[h_{\circledast}^{-1} \Big(u_{\circledast}(a^{h_{\circledast}}_{\circledast, 1}) - \stackrel{\circ}{B^{h}_{\circledast}} \stackrel{\circ}{u^{h}}\Big)
			- \partial_{\circledast}^{h_{\circledast}} u_{\circledcirc}(a^{h_{\circledast}}_{\circledast, 1})\Big]^{2} \mathrm{d}a_{\circledast} \\
			&+ \sum\limits^{}_{\circledast \in \{\MA, \VE\}} \sum_{i_{\circledast} = 2}^{N_{\circledast, h_{\circledast}}}
			\int_{a_{\circledast, i_{\circledast} - 1}^{h_{\circledast}}}^{a_{\circledast, i_{\circledast}}^{h_{\circledast}}}
			\Big(h_{\circledast}^{-1} \big(u_{\circledast}(a^{h_{\circledast}}_{\circledast, i_{\circledast}}) -  u_{\circledast}(a^{h_{\circledast}}_{\circledast, i_{\circledast} - 1})\big)
			- \partial_{\circledast} u_{\circledast}(a^{h_{\circledast}}_{\circledast, i_{\circledast}})\Big)^{2} \mathrm{d}a_{\circledast} \\ 
			&\leq
			2 \sum\limits^{}_{\circledast \in \{\MA, \VE\}} \int_{a_{\circledast, 0}^{h_{\circledast}}}^{a_{\circledast, 1}^{h_{\circledast}}} \Big(u_{\circledast}(0) - \stackrel{\circ}{B^{h}_{\circledast}} \stackrel{\circ}{u^{h}}\Big)^{2} \mathrm{d}a_{\circledast} \\
			&+ 2 \sum\limits^{}_{\circledast \in \{\MA, \VE\}} \sum_{i_{\circledast} = 2}^{N_{\circledast, h_{\circledast}}}
			\int_{a_{\circledast, i_{\circledast} - 1}^{h_{\circledast}}}^{a_{\circledast, i_{\circledast}}^{h_{\circledast}}}
			\Big(h_{\circledast}^{-1} \big(u_{\circledast}(a^{h_{\circledast}}_{\circledast, i_{\circledast}}) -  u_{\circledast}(a^{h_{\circledast}}_{\circledast, i_{\circledast} - 1})\big)
			- \partial_{\circledast} u_{\circledast}(a^{h_{\circledast}}_{\circledast, i_{\circledast}})\Big)^{2} \mathrm{d}a_{\circledast} \to 0
		\end{align*}
		as $h \to 0$.
	\end{enumerate}
	This finishes the proof.
\end{proof}

Let $u^{0} \in D(A)$, $f \in H^{1}(0, T; X) \cap C^{0}\big([0, T], C^{0}(\bar{A}_{\MA}) \times C^{0}(\bar{A}_{\VE})\big)$
and let $u \in C^{1}([0, T], X) \cap C^{0}([0, T], D(A))$ be the corresponding unique classical solution.
Note that we have then
\begin{equation}
	u, \partial_{t} u - Au, f \in C^{0}\big([0, T], C^{0}(\bar{A}_{\MA}) \times C^{0}(\bar{A}_{\VE})\big) \text{ and } u^{0} \in C^{0}(\bar{A}_{\MA}) \times C^{0}(\bar{A}_{\VE}). \notag
\end{equation}
but, in general, not $\partial_{t} u, Au \in C^{0}\big([0, T], C^{0}(\bar{A}_{\MA}) \times C^{0}(\bar{A}_{\VE})\big)$.
Thus, $\partial_{t} u, Au$ cannot be restricted onto the time-space grid
whereas it is possible to restrict $Au$ onto the time grid obtaining an $X$-valued function.

For $\tau > 0$, $h := (h_{\MA}, h_{\VE})$ with $h_{\MA}, h_{\VE} > 0$, we denote
$t \in \stackrel{\circ}{Z^{\tau}}$ and $a_{\circledast} \in \stackrel{\circ}{A_{\circledast}^{h_{\circledast}}}$ for $\circledast \in \{\MA, \VE\}$
\begin{equation}
	u_{\circledast}^{h, \tau}(t, a_{\circledast}) := u_{\circledast}(t, a_{\circledast}), \quad
	f_{\circledast}^{h, \tau}(t, a_{\circledast}) := f_{\circledast}(t, a_{\circledast}), \quad
	u_{\circledast}^{h}(a_{\circledast}) := u^{0}(a_{\circledast}). \notag
\end{equation}
\begin{theorem}[Consistency]
	\label{THEOREM_CONSISTENCY}
	There holds
	\begin{equation}
		\|\mathcal{L}^{h, \tau} u^{h, \tau} - \mathcal{F}(\stackrel{\circ}{u^{0, h}}, \stackrel{\circ}{f^{h, \tau}})\|_{X^{h} \times L^{2}_{\tau}(Z^{\tau}, \stackrel{\circ}{X^{h}}) \times L^{2}_{\tau}(\stackrel{\circ}{Z^{\tau}}, \RR^{2})} \to 0
		\text{ as } (h, \tau) \to 0. \notag
	\end{equation}
\end{theorem}

\begin{proof}
	Splitting the norms of each of the three components, adding and subtracting
	\begin{equation}
		\begin{split}
			\int_{t_{k-1}^{\tau}}^{t_{k}^{\tau}} \int_{a_{\circledast, i_{\circledast} - 1}^{h_{\circledast}}}^{a_{\circledast, i_{\circledast}}^{h_{\circledast}}}
			\Big(&\vartheta \big(\partial_{t} u(t, a_{\circledast}) - Au(t, a_{\circledast} - f(t, a_{\circledast})\big) + \\
			&(1 - \vartheta) \big(\partial_{t} u(t, a_{\circledast}) - Au(t, a_{\circledast}) - f(t, a_{\circledast})\big)\Big) \mathrm{d}a_{\circledast} \mathrm{d}t 
		\end{split}
		\notag
	\end{equation}
	and
	\begin{equation}
		u_{\circledast}(t, 0) - \sum\limits^{}_{\circledcirc \in \{\MA, \VE\}} \int_{A_{\circledcirc}} m_{\circledast \circledcirc}(a_{\circledcirc}) u_{\circledcirc}(t, a_{\circledcirc}) \mathrm{d}a_{\circledcirc} \notag
	\end{equation}
	in the second and third group of terms in
	\begin{equation}
		\|\mathcal{L}^{h, \tau} u^{h, \tau} - \mathcal{F}(\stackrel{\circ}{u^{0, h}}, \stackrel{\circ}{f^{h, \tau}})\|^{2}_{X^{h} \times L^{2}_{\tau}(Z^{\tau}, \stackrel{\circ}{X^{h}}) \times L^{2}_{\tau}(\stackrel{\circ}{Z^{\tau}}, \RR^{2})} \notag
	\end{equation}
	for all $k = 1, \dots, N^{\tau}$, $i_{\circledast} = 1, \dots, N_{\circledast}^{h_{\circledast}}$, $\circledast \in \{\MA, \VE\}$,
	using the definition of $\mathcal{L}^{h, \tau}$ and
	Equations (\ref{EQUATION_MODEL_TRANSFORMED_EQUATION_1})--(\ref{EQUATION_MODEL_TRANSFORMED_EQUATION_6}),
	applying Lemma \ref{LEMMA_CONSISTENCY_FOR_A_H} and Lemma \ref{LEMMA_APPENDIX_DIFFERENCE_OPERATOR_ESTIMATE}
	and exploiting the Cauchy \& Schwarz inequality, we get
	\begin{align*}
		\big\|\mathcal{L}^{h, \tau} &u^{h, \tau} - \mathcal{F}(\stackrel{\circ}{u^{0, h}}, \stackrel{\circ}{f^{h, \tau}})\big\|_{X^{h} \times L^{2}_{\tau}(Z^{\tau}, \stackrel{\circ}{X^{h}}) \times L^{2}_{\tau}(\stackrel{\circ}{Z^{\tau}}, \RR^{2})}^{2} 
		\leq \|u(0, \cdot) - u^{0, h}\|_{\stackrel{\circ}{X^{h}}}^{2} \\
		&\leq 6 \sum_{k = 1}^{N^{\tau}} \int_{t_{k-1}^{\tau}}^{t_{k}^{\tau}} \sum_{\circledast \in \{\MA, \VE\}} \sum_{i_{\circledast} = 1}^{N_{\circledast, h_{\circledast}}} 
		\int_{a_{\circledast, i_{\circledast} - 1}^{h_{\circledast}}}^{a_{\circledast, i_{\circledast}}^{h_{\circledast}}} \big\|\partial_{t} u(t, a_{\circledast}) - A u(t, a_{\circledast}) - f(t, a_{\circledast})\big\|_{\stackrel{\circ}{X^{h}}}^{2} \\
		&+ 6 \sum_{k = 1}^{N^{\tau}} \int_{t_{k-1}^{\tau}}^{t_{k}^{\tau}} \sum_{\circledast \in \{\MA, \VE\}} \sum_{i_{\circledast} = 1}^{N_{\circledast, h_{\circledast}}} 
		\int_{a_{\circledast, i_{\circledast} - 1}^{h_{\circledast}}}^{a_{\circledast, i_{\circledast}}^{h_{\circledast}}}
		\Big(\partial_{t}^{\tau} u(t, a_{\circledast, i_{\circledast}}^{h_{\circledast}}) - \partial_{t} u(t, a_{\circledast}) \Big)^{2} \mathrm{d}a_{\circledast} \mathrm{d}t \\
		&+ 6 \vartheta \sum_{k = 1}^{N^{\tau}} \int_{t_{k-1}^{\tau}}^{t_{k}^{\tau}} \sum\limits^{}_{\circledast \in \{\MA, \VE\}} \sum\limits_{i_{\circledast} = 1}^{N_{\circledast, h_{\circledast}}}
		\int\limits_{a_{\circledast, i_{\circledast} - 1}^{h_{\circledast}}}^{a_{\circledast, i_{\circledast}}^{h_{\circledast}}}
		\Big(A_{\circledast} u(t, a_{\circledast}) - \stackrel{\circ}{A_{\circledast}} \stackrel{\circ}{u}_{\circledast}(t_{k}^{\tau}, a_{\circledast, i_{\circledast}}^{h_{\circledast}})\Big)^{2} \mathrm{d}a_{\circledast} \mathrm{d}t \\
		&+ 6 (1 - \vartheta) \sum_{k = 1}^{N^{\tau}} \int_{t_{k-1}^{\tau}}^{t_{k}^{\tau}} \sum\limits^{}_{\circledast \in \{\MA, \VE\}} \sum\limits_{i_{\circledast} = 1}^{N_{\circledast, h_{\circledast}}}
		\int\limits_{a_{\circledast, i_{\circledast} - 1}^{h_{\circledast}}}^{a_{\circledast, i_{\circledast}}^{h_{\circledast}}}
		\Big(A_{\circledast} u(t, a_{\circledast}) - \stackrel{\circ}{A_{\circledast}} \stackrel{\circ}{u}_{\circledast}(t_{k-1}^{\tau}, a_{\circledast, i_{\circledast}}^{h_{\circledast}})\Big)^{2} \mathrm{d}a_{\circledast} \mathrm{d}t \displaybreak \\
		&+ 6 \vartheta \sum_{k = 1}^{N^{\tau}} \int_{t_{k-1}^{\tau}}^{t_{k}^{\tau}} \sum_{\circledast \in \{\MA, \VE\}} \sum_{i_{\circledast} = 1}^{N_{\circledast, h_{\circledast}}}
		\int_{a_{\circledast, i_{\circledast} - 1}^{h_{\circledast}}}^{a_{\circledast, i_{\circledast}}^{h_{\circledast}}}
		\Big(\stackrel{\circ}{f^{h, \tau}}(t_{k}^{\tau}, a_{\circledast, i_{\circledast}}^{h_{\circledast}}) - f(t, a_{\circledast})\Big)^{2} \mathrm{d}a_{\circledast} \mathrm{d}t \\
		&+ 6 (1 - \vartheta) \sum_{k = 1}^{N^{\tau}} \int_{t_{k-1}^{\tau}}^{t_{k}^{\tau}} \sum_{\circledast \in \{\MA, \VE\}} \sum_{i_{\circledast} = 1}^{N_{\circledast, h_{\circledast}}} 
		\int_{a_{\circledast, i_{\circledast} - 1}^{h_{\circledast}}}^{a_{\circledast, i_{\circledast}}^{h_{\circledast}}}
		\Big(\stackrel{\circ}{f^{h, \tau}}(t_{k-1}^{\tau}, a_{\circledast, i_{\circledast}}^{h_{\circledast}}) - f(t, a_{\circledast})\Big)^{2} \mathrm{d}a_{\circledast} \mathrm{d}t \\
		&+ 2 \sum_{k = 1}^{N^{\tau}} \int_{t_{k-1}^{\tau}}^{t_{k}^{\tau}} \sum_{\circledast \in \{\MA, \VE\}}
		\Big(\stackrel{\circ} {B_{\circledast}^{h}} \stackrel{\circ}{u^{h}}(t, a_{\circledast}) - \sum_{\circledcirc \in \{\MA, \VE\}} \int_{0}^{a_{\ast}^{\dag}} m_{\circledast \circledcirc}(a_{\circledcirc}) \mathrm{d} a_{\circledcirc}\Big)^{2} \mathrm{d}t \\
		&+ 2 \sum_{k = 1}^{N^{\tau}} \int_{t_{k-1}^{\tau}}^{t_{k}^{\tau}} \sum_{\circledast \in \{\MA, \VE\}}
		\Big(u_{\circledast}(t, 0) - \sum\limits^{}_{\circledcirc \in \{\MA, \VE\}} \int_{A_{\circledcirc}} m_{\circledast \circledcirc}(a_{\circledcirc}) u_{\circledcirc}(t, a_{\circledcirc}) \mathrm{d}a_{\circledcirc}\Big)^{2} \mathrm{d}a_{\circledcirc} \mathrm{d}t
		\to 0
	\end{align*}
	as $(h, \tau) \to 0$.
\end{proof}

\subsection{Stability and Convergence}
Our stability investigations are very much related to deducing a resolvent estimate in Section \ref{SECTION_WELL_POSEDNESS_AND_ASYMPTOTICS}.
Whereas the latter was obtained using multiplier techniques based on partial integration,
a summation by parts formula will be expoloited here to obtain a uniform resolvent estimate for $\stackrel{\circ}{A^{h}}$.
Further, a uniform $L^{\infty}$-estimate for the numerical solution based on the rational approximation 
for the corresponding $C_{0}$-semigroup will be shown.
Together with the consistency result from the previous subsection,
this will lead to the unconditional convergence of the implicit scheme.

We let
\begin{equation}
	\omega_{0} := \max\{a_{\MA}^{\dag}, a_{\VE}^{\dag}, \tfrac{1}{2} M^{2}\} > 0 \notag
\end{equation}
for $M :=
\max\limits^{}_{\circledast, \circledcirc \in \{\MA, \VE\}}
\max\limits^{}_{a_{\circledcirc} \in A_{\circledcirc}^{h_{\circledcirc}}} |m_{\circledast \circledcirc}(a_{\circledcirc})| \leq
\max\limits^{}_{\circledast, \circledcirc \in \{\MA, \VE\}} \|m_{\circledast \circledcirc}\|_{L^{\infty}(A_{\circledast})} < \infty$.\\[-0.5cm]
\begin{lemma}
	\label{LEMMA_ESTIMATE_FOR_A_H}
	For any $h_{\MA}, h_{\VE} > 0$, there holds for any $u^{h} \in D(\stackrel{\circ}{A^{h}})$
	\begin{equation}
		\langle \stackrel{\circ}{A^{h}} u^{h}, u^{h}\rangle_{\stackrel{\circ}{X^{h}}} \leq \omega_{0} \|u^{h}\|_{\stackrel{\circ}{X^{h}}}^{2}. \notag
	\end{equation}
\end{lemma}

\begin{proof}
	Let $u^{h} \in D(\stackrel{\circ}{A^{h}})$ and let $\beta \geq 0$ be an arbitrary number to be fixed latter.
	Using Lemma \ref{LEMMA_SUMMATION_BY_PARTS}, we can estimate
	\begin{equation}
		-\sum_{i_{\circledast} = 2}^{N_{\circledast}^{h_{\circledast}}} \big(\partial_{\circledast}^{h_{\circledast}} u^{h_{\circledast}}\big)_{\circledast, i_{\circledast}} u_{\circledast, i_{\circledast}}^{h_{\circledast}}
		\leq \tfrac{1}{2} h_{\circledast}^{-1}
		\Big(\big(u_{\circledast, 1}^{h_{\circledast}}\big)^{2} - \big(u_{\circledast, N_{\circledast}^{h_{\circledast}}}^{h_{\circledast}}\big)^{2}\Big).
		\notag
	\end{equation}
	Hence,
	\begin{align*}
		\langle \stackrel{\circ}{A^{h}} u^{h} - &\beta u^{h}, u^{h}\rangle_{\stackrel{\circ}{X^{h}}} = \langle \stackrel{\circ}{A^{h}} u^{h}, u^{h}\rangle_{\stackrel{\circ}{X^{h}}} - \beta \|u^{h}\|_{\stackrel{\circ}{X^{h}}}^{2} \\
		&= -\sum_{\circledast \in \{\MA, \VE\}} \Big(u_{\circledast, 1}^{h_{\circledast}} -
		\sum_{\circledcirc \in \{\MA, \VE\}} h_{\circledcirc} \sum\limits_{i_{\circledcirc} = 1}^{N_{\circledcirc, h_{\circledast}}} m_{\circledast \circledcirc}(a_{\circledcirc, i_{\circledcirc}}^{h_{\circledcirc}}) u_{\circledcirc, i_{\circledcirc}}^{h_{\circledcirc}}\Big) u_{\circledast, 1}^{h_{\circledast}} \displaybreak \\
		&\phantom{=}\; -\sum_{\circledast \in \{\MA, \VE\}} h_{\circledast} \sum_{i_{\circledast} = 2}^{N_{\circledast}^{h_{\circledast}}} \big(\partial_{a_{\circledast}}^{h_{\circledast}} u^{h}_{\circledast, i_{\circledast}}\big) u^{h}_{\circledast, i_{\circledast}}
		- \beta \|u^{h}\|_{\stackrel{\circ}{X^{h}}}^{2} \\
		&\leq \sum_{\circledast \in \{\MA, \VE\}}
		\Big(-\tfrac{1}{2} \big(u_{\circledast, 1}^{h_{\circledast}}\big)^{2} + \sum_{\circledcirc \in \{\MA, \VE\}} h_{\circledcirc} \sum\limits_{i_{\circledcirc} = 1}^{N_{\circledcirc, h_{\circledast}}} m_{\circledast \circledcirc}(a_{\circledcirc, i_{\circledcirc}}^{h_{\circledcirc}}) u_{\circledcirc, i_{\circledcirc}}^{h_{\circledcirc}} u_{\circledast, 1}^{h_{\circledast}}\Big)
		- \beta \|u^{h}\|_{\stackrel{\circ}{X^{h}}}^{2} \\
		&\leq \tfrac{1}{2} M^{2} \sum_{\circledast, \circledcirc \in \{\MA, \VE\}} h_{\circledcirc} \sum\limits_{i_{\circledcirc} = 1}^{N_{\circledcirc, h_{\circledast}}} \big(u_{\circledcirc, i_{\circledcirc}}^{h_{\circledcirc}}\big)^{2}
		- \beta \|u^{h}\|_{\stackrel{\circ}{X^{h}}}^{2} \\
		&\leq (\tfrac{1}{2} M^{2} - \beta) \|u^{h}\|_{\stackrel{\circ}{X^{h}}}^{2}.
	\end{align*}
	The claim follows now for $\beta := \omega_{0}$.
\end{proof}

From Lemma \ref{LEMMA_ESTIMATE_FOR_A_H}, we get using \cite[Theorem 4.2]{Pa1983} the following resolvent estimate for $\stackrel{\circ}{A^{h}}$.
\begin{corollary}
	\label{COROLLARY_ESTIMATE_FOR_THE_RESOLVENT_OF_A_H}
	For $\lambda \in (\omega_{0}, \infty)$, the operator $\lambda \mathrm{id} - \stackrel{\circ}{A^{h}}$ is continuously invertible with
	\begin{equation}
		\|(\lambda \mathrm{id} - \stackrel{\circ}{A^{h}})^{-1}\|_{L(\stackrel{\circ}{X^{h}})} \leq (\lambda - \omega_{0})^{-1}. \notag
	\end{equation}
\end{corollary}

Now, we can prove the following unconditional stability result.
\begin{theorem}[Stability]
	\label{THEOREM_STABILITY}
	Let $\vartheta \in [\tfrac{1}{2}, 1]$ and let $\bar{\tau} := \tfrac{1}{2 \vartheta \omega_{0}}$.
	For any $\tau \in (0, \bar{\tau})$ and $h = (h_{\MA}, h_{\VE})$ with $h_{\MA}, h_{\VE} > 0$, there exists an number $C > 0$ such that
	any data $u^{0, h} \in X^{h}$, $\stackrel{\circ}{f^{h, \tau}} \in L^{2}_{\tau}(Z^{\tau}, \stackrel{\circ}{X^{h}})$ admit a unique numerical solution
	$u^{h, \tau} \in H^{1}_{\tau}\big(Z^{\tau}, H^{1}_{h_{\MA}}(A_{\MA}^{h_{\MA}}) \times H^{1}_{h_{\VE}}(A_{\VE}^{h_{\VE}})\big)$ to Equation (\ref{EQUATION_NUMERICAL_SCHEME_ABSTRACT_FORM})
	depending continuously on the data in terms of the estimate
	\begin{equation}
		\|u^{h, \tau}\|_{L^{\infty}_{\tau}(Z^{\tau}, X^{h})}
		\leq C \big\|\mathcal{F} \big(u^{0, h}, f^{h, \tau}\big)\big\|_{\stackrel{\circ}{X^{h}} \times L^{2}_{\tau}(\stackrel{\circ}{Z^{\tau}}, \stackrel{\circ}{X^{h}}) \times L^{2}_{\tau}(\stackrel{\circ}{Z^{\tau}}, \RR^{2})}.
		\notag
	\end{equation}
\end{theorem}

\begin{proof}
	Recalling that Equations (\ref{EQUATION_NUMERICAL_SCHEME_ABSTRACT_FORM}) and
	(\ref{EQUATION_MODEL_TRANSFORMED_DISCRETIZED_EQUATION_1})--(\ref{EQUATION_MODEL_TRANSFORMED_DISCRETIZED_EQUATION_3}) are equivalent,
	Equation (\ref{EQUATION_NUMERICAL_SCHEME_ABSTRACT_FORM}) can be written as
	\begin{align}
		u^{h, \tau}(0, \cdot) &= u^{0, h}, \label{EQUATION_FULL_DISCRETIZATION_ITERATION_1} \\
		\big(\tfrac{1}{\tau} - \vartheta \stackrel{\circ}{A^{h}}\big) \stackrel{\circ}{u^{h, \tau}}(t_{k}, \cdot) &=
		\big(\tfrac{1}{\tau} + (1 - \vartheta) \stackrel{\circ}{A^{h}}\big) \stackrel{\circ}{u^{h, \tau}}(t_{k-1}, \cdot) \label{EQUATION_FULL_DISCRETIZATION_ITERATION_2} \\
		&+ \vartheta \stackrel{\circ}{f^{h, \tau}}(t_{k}, \cdot) + (1 - \vartheta) \stackrel{\circ}{f^{h, \tau}}(t_{k-1}, \cdot) \text{ for } k = 1, \dots, N^{\tau}, \notag \\
		\stackrel{\diamond}{u^{h, \tau}}(t_{k}, \cdot) &= \;\stackrel{\circ}{B^{h}} \stackrel{\circ}{u^{h, \tau}}(t_{k}, \cdot) \text{ for } k = 1, \dots, N^{\tau}. \label{EQUATION_FULL_DISCRETIZATION_ITERATION_3}
	\end{align}
	One can easily observe that Equations (\ref{EQUATION_FULL_DISCRETIZATION_ITERATION_1})--(\ref{EQUATION_FULL_DISCRETIZATION_ITERATION_2})
	and (\ref{EQUATION_FULL_DISCRETIZATION_ITERATION_3}) decouple.
	Given a solution to the difference equations (\ref{EQUATION_FULL_DISCRETIZATION_ITERATION_1})--(\ref{EQUATION_FULL_DISCRETIZATION_ITERATION_2}),
	a solution to Equation (\ref{EQUATION_FULL_DISCRETIZATION_ITERATION_3}) can explicitly obtained.
	Thus, Equations (\ref{EQUATION_MODEL_TRANSFORMED_DISCRETIZED_EQUATION_1})--(\ref{EQUATION_MODEL_TRANSFORMED_DISCRETIZED_EQUATION_3})
	are uniquely solvable if this is the case for Equations (\ref{EQUATION_FULL_DISCRETIZATION_ITERATION_1})--(\ref{EQUATION_FULL_DISCRETIZATION_ITERATION_2}).
	The latter are uniquely solvable for any data if and only if the operator $\mathrm{id} - \tau \vartheta \stackrel{\circ}{A^{h}}$ is (continuously) invertible.
	According to Corollary \ref{COROLLARY_ESTIMATE_FOR_THE_RESOLVENT_OF_A_H},
	the latter is the case if $\tau \in (0, \tfrac{\omega_{0}}{\vartheta})$.

	Letting
	\begin{equation}
		H_{1}^{h, \tau} := \big(\tfrac{1}{\tau} \mathrm{id} - \vartheta \stackrel{\circ}{A^{h}}\big), \quad
		H_{2}^{h, \tau} := \big(\tfrac{1}{\tau} \mathrm{id} + (1 - \vartheta) \stackrel{\circ}{A^{h}}\big), \notag
	\end{equation}
	we can easily show by induction that the unique solution to Equation (\ref{EQUATION_FULL_DISCRETIZATION_ITERATION_1})--(\ref{EQUATION_FULL_DISCRETIZATION_ITERATION_2})
	is iteratively given by
	\begin{equation}
		\begin{split}
			\stackrel{\circ}{u^{h, \tau}}(t_{k}^{\tau}, \cdot) &=
			\big((H_{1}^{h, \tau})^{-1} H_{2}^{h, \tau}\big)^{k} \stackrel{\circ}{u^{0, h}} \\
			&+ \sum_{j = 1}^{k}
			\big((H_{1}^{h, \tau})^{-1} H_{2}^{h, \tau}\big)^{k - j} (H_{1}^{h, \tau})^{-1}
			\big(\vartheta f^{\tau, h}(t_{k}^{\tau}, \cdot) + (1 - \vartheta) \stackrel{\circ}{f^{\tau, h}}(t_{k-1}^{\tau}, \cdot)\big) \\
			&\text{ for } k = 0, \dots, N^{\tau}.
		\end{split}
		\label{EQUATION_OPERATOR_EQUATION_FOR_L_H_TAU_EXPLICIT_FORM_SOLUTION}
	\end{equation}
	Further, for $\tau \in (0, \tfrac{\omega_{0}}{\vartheta})$, we trivially obtain the operator identity
	\begin{equation}
		\begin{split}
			(H_{1}^{h, \tau})^{-1} H_{2}^{h, \tau} &=
			\big(\tfrac{1}{\tau} \mathrm{id} - \vartheta \stackrel{\circ}{A^{h}}\big)^{-1} \big(\tfrac{1}{\tau} \mathrm{id} + (1 - \vartheta) \stackrel{\circ}{A^{h}}\big) \\
			&= \mathrm{id} - \big(\tfrac{1}{\tau} \mathrm{id} - \vartheta \stackrel{\circ}{A^{h}}\big)^{-1} \stackrel{\circ}{A^{h}} \\
			&= \mathrm{id} - \tfrac{1}{\vartheta} \big(\tfrac{1}{\tau} \mathrm{id} - \vartheta \stackrel{\circ}{A^{h}}\big)^{-1} \big(\tfrac{1}{\tau} \mathrm{id} - \vartheta \stackrel{\circ}{A^{h}} - \tfrac{1}{\tau} \mathrm{id}\big) \\
			&= \mathrm{id} - \tfrac{1}{\vartheta} \Big(\mathrm{id} - \tfrac{1}{\tau} \big(\tfrac{1}{\tau} - \vartheta \stackrel{\circ}{A^{h}}\big)^{-1}\Big) \\
			&= \mathrm{id} - \tfrac{1}{\vartheta} \Big(\mathrm{id} - \big(\mathrm{id} - \tau \vartheta \stackrel{\circ}{A^{h}}\big)^{-1}\Big) \\
			&= (1 - \tfrac{1}{\vartheta}) \mathrm{id} + \big(\mathrm{id} - \tau \vartheta \stackrel{\circ}{A^{h}}\big)^{-1}. 
		\end{split}
		\label{EQUATION_ESTIMATE_OPERATOR_H_1_TIMES_H_2}
	\end{equation}
	Using again Corollary \ref{COROLLARY_ESTIMATE_FOR_THE_RESOLVENT_OF_A_H}, we can estimate for $\tau \in (0, \tfrac{1}{2 \vartheta \omega_{0}})$
	\begin{equation}
		\|(H_{1}^{h, \tau})^{-1}\|_{L(\stackrel{\circ}{X^{h}})} = \|(\tfrac{1}{\tau} \mathrm{id} - \vartheta \stackrel{\circ}{A^{h}})^{-1}\|_{L(\stackrel{\circ}{X^{h}})} \leq
		\tfrac{1}{\vartheta} \|(\tfrac{1}{\vartheta \tau} - \stackrel{\circ}{A^{h}})^{-1}\|_{L(\stackrel{\circ}{X^{h}})} \leq 4 \tau. \notag
	\end{equation}
	This together with Equation (\ref{EQUATION_ESTIMATE_OPERATOR_H_1_TIMES_H_2}) implies
	\begin{equation}
		\|(H_{1}^{h, \tau})^{-1} H_{2}^{h, \tau}\|_{L(\stackrel{\circ}{X^{h}})}
		\leq 1 + 4 \tau \text{ for } \tau \in (0, \tfrac{1}{2 \vartheta \omega_{0}}). \notag
	\end{equation}
	Therefore, for any $k = 0, \dots, N^{\tau}$,
	\begin{equation}
		\big\|\big((H_{1}^{h, \tau})^{-1} H_{2}^{h, \tau}\big)^{k}\big\|_{L(\stackrel{\circ}{X^{h}})}
		\leq \sqrt{2} (1 + 4 \tau)^{k} \leq \exp(4 \tau k) \leq \exp(4T). \notag
	\end{equation}
	Recalling now Equation (\ref{EQUATION_OPERATOR_EQUATION_FOR_L_H_TAU_EXPLICIT_FORM_SOLUTION}) and
	applying Young's inequality, we obtain for all $k = 0, \dots, N^{\tau}$
	\begin{equation}
		\begin{split}
			\|\stackrel{\circ}{u^{h, \tau}}(t_{k}^{\tau}, \cdot)\|_{L(\stackrel{\circ}{X^{h}})} &\leq
			\exp(4T) \Big(\|\stackrel{\circ}{u^{0, h}}\|_{L(\stackrel{\circ}{X^{h}})} \\
			&+ 4\tau \sum_{j = 1}^{k} \big\|\big(\vartheta \stackrel{\circ}{f^{\tau, h}}(t_{k-j}^{\tau}, \cdot)
			+ (1 - \vartheta) \stackrel{\circ}{f^{\tau, h}}(t_{k-j-1}^{\tau}, \cdot)\big)\big\|_{L(\stackrel{\circ}{X^{h}})}\Big) \\
			&\leq \sqrt{2} (1 + 4\tau) \exp(4T) \big\|\mathcal{F} \big(u^{0, h}, f^{h, \tau}\big)\big\|_{\stackrel{\circ}{X^{h}} \times L^{2}_{\tau}(\stackrel{\circ}{Z^{\tau}}, \stackrel{\circ}{X^{h}}) \times L^{2}_{\tau}(\stackrel{\circ}{Z^{\tau}}, \RR^{2})}.
		\end{split}
		\label{EQUATION_STABILITY_ESTIMATE_1}
	\end{equation}
	Next, Equation (\ref{EQUATION_FULL_DISCRETIZATION_ITERATION_3}) uniquely determines the unknown $\stackrel{\diamond}{u^{h, \tau}}$
	which we can easily be estimated as follows
	\begin{equation}
		\begin{split}
			\|\stackrel{\diamond}{u^{h, \tau}}(t_{k}^{\tau})\|_{\RR^{2}} &\leq
			\|\stackrel{\circ}{B^{h}}\|_{L(\stackrel{\circ}{X^{h}})} \|\stackrel{\circ}{u^{h}}(t_{k}^{\tau}, \cdot)\|_{L(\stackrel{\circ}{X^{h}})} \\
			&\leq \sqrt{2} M (1 + 4\tau) \exp(4T) \big\|\mathcal{F} \big(u^{0, h}, f^{h, \tau}\big)\big\|_{\stackrel{\circ}{X^{h}} \times L^{2}_{\tau}(\stackrel{\circ}{Z^{\tau}}, \stackrel{\circ}{X^{h}}) \times L^{2}_{\tau}(\stackrel{\circ}{Z^{\tau}}, \RR^{2})}.
		\end{split}
		\label{EQUATION_STABILITY_ESTIMATE_2}
	\end{equation}
	Estimates from Equations (\ref{EQUATION_STABILITY_ESTIMATE_1}), (\ref{EQUATION_STABILITY_ESTIMATE_2})
	together with Young's inequality yield now the claim
	with $C := 2 (1 + M) (1 + 4\tau) \exp(4T)$.
\end{proof}

Again, let $u^{0} \in D(A)$, $f \in H^{1}(0, T; X) \cap C^{0}\big([0, T], C^{0}(\bar{A}_{\MA}) \times C^{0}(\bar{A}_{\VE})\big)$
and let $\bar{u} \in H^{1}(0, T; X) \cap L^{2}(0, T; D(A))$ be the corresponding unique classical solution.
For $\vartheta \in [\tfrac{1}{2}, 1]$, $\tau > 0$ and $h := (h_{\MA}, h_{\VE})$ satisfying the conditions of Theorem \ref{THEOREM_STABILITY},
$t \in \stackrel{\circ}{Z^{\tau}}$ and $a_{\circledast} \in \stackrel{\circ}{A_{\circledast}^{h_{\circledast}}}$ for $\circledast \in \{\MA, \VE\}$
\begin{equation}
	f_{\circledast}^{h, \tau}(t, a_{\circledast}) := f_{\circledast}(t, a_{\circledast}), \quad
	u_{\circledast}^{h}(a_{\circledast}) := u^{0}(a_{\circledast}). \notag
\end{equation}
Let $u^{h, \tau} \in H^{1}\big(Z^{\tau}, H^{1}_{h_{\MA}}(A_{\MA}^{h_{\MA}}) \times H^{1}_{h_{\VE}}(A_{\VE}^{h_{\VE}})\big)$ denote the unique solution of Equation (\ref{EQUATION_NUMERICAL_SCHEME_ABSTRACT_FORM})
given in Theorem \ref{THEOREM_STABILITY}.

Using the Lax' principle, we have
\begin{theorem}[Convergence]
	There holds
	\begin{equation}
		\big\|\bar{u}^{h, \tau} - u^{h, \tau}\big\|_{L^{\infty}_{\tau}(Z^{\tau}, X^{h})} \to 0 \text{ as } (h, \tau) \to 0. \notag
	\end{equation}
\end{theorem}

\begin{proof}
	Using the fact that $\mathcal{L}^{h, \tau} \bar{u}^{h, \tau} = \mathcal{F}^{h, \tau}\big(\stackrel{\circ}{u^{0, h}}, \stackrel{\circ}{f^{h, \tau}}\big)$
	and exploiting Theorems \ref{THEOREM_STABILITY} and \ref{THEOREM_CONSISTENCY}, we get
	\begin{align*}
		\big\|\bar{u}^{h, \tau} - u^{h, \tau}\big\|_{L^{\infty}_{\tau}(Z^{\tau}, X^{h})} &\leq
		C \big\|\mathcal{L}^{h, \tau} \big(\bar{u}^{h, \tau} - u^{h, \tau}\big)\big\|_{\stackrel{\circ}{X^{h}} \times L^{2}_{\tau}(\stackrel{\circ}{Z^{\tau}}, \stackrel{\circ}{X^{h}}) \times L^{2}_{\tau}(\stackrel{\circ}{Z^{\tau}}, \RR^{2})} \\
		&= C \big\|\mathcal{L}^{h, \tau} \bar{u}^{h, \tau}\big\|_{\stackrel{\circ}{X^{h}} \times L^{2}_{\tau}(\stackrel{\circ}{Z^{\tau}}, \stackrel{\circ}{X^{h}}) \times L^{2}_{\tau}(\stackrel{\circ}{Z^{\tau}}, \RR^{2})} \to 0 \notag
	\end{align*}
	as $(h, \tau) \to 0$.
\end{proof}

\subsection{Computer Implementation and Numerical Example}
\label{IMPLEMENATION_AND_EXAMPLE}
In this Section, we use our developments from the previous Section \ref{SECTION_FINITE_DIFFERENCE_SCHEME}
and construct an algorithm to numerically solve Equations (\ref{EQUATION_MODEL_TRANSFORMED_EQUATION_1})--(\ref{EQUATION_MODEL_TRANSFORMED_EQUATION_6}).
Throughout this Section, all discrete spaces will be viewed as the usual Euclidian ones
and all linear operators will be replaced by matrices.
In particular,
\begin{equation}
	\stackrel{\circ}{X^{h}} \simeq \RR^{N_{\MA}} \times \RR^{N_{\VE}} \simeq \RR^{N_{\MA}^{h_{\MA}} + N_{\VE}^{h_{\VE}}} \text{ and }
	L^{2}_{\tau}(Z^{\tau}, \stackrel{\circ}{X^{h}}) \simeq \RR^{N^{\tau} \times (N_{\MA}^{h_{\MA}} + N_{\VE}^{h_{\VE}})}, \text{ etc.}
	\notag
\end{equation}
Introducing the matrices
\begin{equation}
	\stackrel{\circ}{B_{\circledast \circledcirc}^{h}} := {\footnotesize
	\begin{pmatrix}
		m_{\circledast \circledcirc}(a_{\circledcirc, 1}^{h_{\circledcirc}}) & m_{\circledast \circledcirc}(a_{\circledcirc, 2}^{h_{\circledcirc}}) & m_{\circledast \circledcirc}(a_{\circledcirc, 3}^{h_{\circledcirc}}) & \hdots & m_{\circledast \circledcirc}(a_{\circledcirc, N_{\circledcirc}^{h_{\circledcirc}}}^{h_{\circledcirc}}) \\
	\end{pmatrix}} \text{ for } \circledast, \circledcirc \in \{\MA, \VE\} \notag
\end{equation}
and
\begin{align*}
	\stackrel{\circ}{A_{\MA \MA}^{h}} &:= {\footnotesize
	\begin{pmatrix}
		-\tfrac{1}{h_{\MA}} + m_{\MA \MA}(a_{\MA, 1}^{h_{\MA}}) & m_{\MA \MA}(a_{\MA, 2}^{h_{\MA}}) & m_{\MA \MA}(a_{\MA, 3}^{h_{\MA}}) & \hdots & m_{\MA \MA}(a_{\MA, N_{\MA}^{h_{\MA}} - 1}^{h_{\MA}}) & m_{\MA \MA}(a_{\MA, N_{\MA}^{h_{\MA}}}^{h_{\MA}}) \\
		-\tfrac{1}{h_{\MA}} & \tfrac{1}{h_{\MA}} & 0 & \hdots & 0 & 0  \\
		0 & -\tfrac{1}{h_{\MA}} & \tfrac{1}{h_{\MA}} & \hdots & 0 & 0 \\
		\vdots & \vdots & \vdots & \ddots & \vdots & \vdots \\
		0 & 0 & 0 & \hdots & -\tfrac{1}{h_{\MA}} & \tfrac{1}{h_{\MA}}
	\end{pmatrix}}, \\
	\stackrel{\circ}{A_{\VE \VE}^{h}} &:= {\footnotesize
	\begin{pmatrix}
		-\tfrac{1}{h_{\VE}} + m_{\VE \VE}(a_{\VE, 1}^{h_{\VE}}) & m_{\VE \VE}(a_{\VE, 2}^{h_{\VE}}) & m_{\VE \VE}(a_{\VE, 3}^{h_{\VE}}) & \hdots & m_{\VE \VE}(a_{\VE, N_{\VE}^{h_{\VE}} - 1}^{h_{\VE}}) & m_{\VE \VE}(a_{\VE, N_{\VE}^{h_{\VE}}}^{h_{\VE}}) \\
		-\tfrac{1}{h_{\VE}} & \tfrac{1}{h_{\VE}} & 0 & \hdots & 0 & 0  \\
		0 & -\tfrac{1}{h_{\VE}} & \tfrac{1}{h_{\VE}} & \hdots & 0 & 0 \\
		\vdots & \vdots & \vdots & \ddots & \vdots & \vdots \\
		0 & 0 & 0 & \hdots & -\tfrac{1}{h_{\VE}} & \tfrac{1}{h_{\VE}}
	\end{pmatrix}}, \\
	\stackrel{\circ}{A_{\circledast \circledcirc}^{h}} &= {\footnotesize
	\frac{h_{\circledcirc}}{h_{\circledast}}
	\begin{pmatrix}
		m_{\circledast \circledcirc}(a_{\circledcirc, 1}^{h_{\circledcirc}}) & m_{\circledast \circledcirc}(a_{\circledcirc, 2}^{h_{\circledcirc}}) & m_{\circledast \circledcirc}(a_{\circledcirc, 3}^{h_{\circledcirc}}) & \hdots & m_{\circledast \circledcirc}(a_{\circledcirc, N_{\circledcirc}^{h_{\circledcirc}}}^{h_{\circledcirc}}) \\
		0 & 0 & 0 & \hdots & 0 \\
		\vdots & \vdots & \vdots & \ddots & \vdots \\
		0 & 0 & 0 & \hdots & 0
	\end{pmatrix}} \text{ for } \{\circledast, \circledcirc\} = \{\MA, \VE\},
\end{align*}
the operators $\stackrel{\circ}{A^{h}}$ and $\stackrel{\circ}{B^{h}}$ can be represented in the matrix form
\begin{equation}
	\stackrel{\circ}{A^{h}} =
	\begin{pmatrix}
		\stackrel{\circ}{A}_{\MA \MA}^{h} & \stackrel{\circ}{A}_{\MA \VE}^{h} \\
		\stackrel{\circ}{A}_{\VE \MA}^{h} & \stackrel{\circ}{A}_{\VE \VE}^{h}
	\end{pmatrix} \in \RR^{N \times N} \text{ and }
	\stackrel{\circ}{B^{h}} =
	\begin{pmatrix}
		\stackrel{\circ}{B}_{\MA \MA}^{h} & \stackrel{\circ}{B}_{\MA \VE}^{h} \\
		\stackrel{\circ}{B}_{\VE \MA}^{h} & \stackrel{\circ}{B}_{\VE \VE}^{h}
	\end{pmatrix} \in \RR^{2 \times N}
	\text{ with } N := N_{\MA}^{h_{\MA}} + N_{\VE}^{h_{\VE}}. \notag
\end{equation}
Further, we write
$u^{h, \tau; k}$, $\stackrel{\circ}{u^{h, \tau; k}}$, $\stackrel{\diamond}{u^{h, \tau; k}}$ and $\stackrel{\circ}{f^{h, \tau; k}}$
for $u^{h, \tau}(t_{k}, \cdot)$, $\stackrel{\circ}{u^{h, \tau}}(t_{k}, \cdot)$, $\stackrel{\diamond}{u^{h, \tau}}(t_{k}, \cdot)$ and $\stackrel{\circ}{f^{h, \tau}}(t_{k}, \cdot)$, respectively,
for $k = 0, \dots, N^{\tau}$.
With this notation,
Equations (\ref{EQUATION_FULL_DISCRETIZATION_1})--(\ref{EQUATION_FULL_DISCRETIZATION_3})
reduce to a system of linear algebraic equations
\begin{align}
	\Big(\tfrac{1}{\tau} - \vartheta \stackrel{\circ}{A^{h}}\Big) {\stackrel{\circ}{u}}^{h, \tau; k} &=
	\Big(\tfrac{1}{\tau} + (1 - \vartheta) \stackrel{\circ}{A^{h}}\Big) {\stackrel{\circ}{u}}^{h, \tau; k-1} \label{EQUATION_FULL_DISCRETIZATION_ALGEBRAIC_1} \\
	&+ \vartheta \tau {\stackrel{\circ}{f}}^{h, \tau; k} + (1 - \vartheta) \tau {\stackrel{\circ}{f}}^{h, \tau; k-1} \text{ for } k = 1, \dots, N^{\tau}, \notag \\
	{\stackrel{\diamond}{u}}^{h, \tau; k} &= \; \stackrel{\circ}{B^{h}} {\stackrel{\circ}{u}}^{h, \tau; k-1} \text{ for } k = 1, \dots, N^{\tau}, \label{EQUATION_FULL_DISCRETIZATION_ALGEBRAIC_2} \\
	{\stackrel{\circ}{u}}^{h, \tau; 0} &= \; \stackrel{\circ}{u^{0, h}}, \quad 
	{\stackrel{\diamond}{u}}^{h, \tau; 0} = \; \stackrel{\diamond}{u^{0, h}}. \label{EQUATION_FULL_DISCRETIZATION_ALGEBRAIC_3}
\end{align}
By the virtue of Theorem \ref{THEOREM_STABILITY},
Equation (\ref{EQUATION_FULL_DISCRETIZATION_ALGEBRAIC_1})--(\ref{EQUATION_FULL_DISCRETIZATION_ALGEBRAIC_3}) is uniquely solvable
if $\vartheta \in [1/2, 1]$ and $\tau \in (0, \tfrac{1}{2 \vartheta \omega_{0}})$.

\subsection{U.S. Population in 2011: Reported vs. Simulated}
To verify our model and test the numerical scheme,
we ran a numerical simulation
to predict the growth of the United States population over the decade between 2001 and 2011.
The information on the population structure in 2001 and 2011
was obtained from the International Data Base
of the U.S. Bureau of Census \cite{USCB2013} (last updated in December 2013).

During the whole period of 2001--2011,
the age-specific survival probabilities both for men and women were assumed to be constantly equal to those reported for 2011 in \cite[Table 1, pp. 202--203]{Chu2007}.
The birth rates by age of mother were selected to be constantly equal to those reported for 2008 in \cite[Table 4, p. 52]{Ma2013}.
The sex ratio was chosen as 1.05 (cf. \cite{CIA2014}).
The annual net immigration was selected as the average net immigration over the period 2001--2009 as reported in \cite[Table 2]{ShrHei2011}.
Due to the lack of more accurate information,
the age and sex structure of the newcomer immigrants' cohort was assumed to be the same
as of those immigrants who have already dwelled in the U.S. in 2001 or before (see \cite{MPI2014}).
Unless the data were divided into single-year age groups,
the average value in each of the groups was computed to estimate each of the single-year values.

Using the age-specific survival probabilities,
all system data and parameters were transformed to the form (\ref{EQUATION_FULL_DISCRETIZATION_ALGEBRAIC_1})--(\ref{EQUATION_FULL_DISCRETIZATION_ALGEBRAIC_3}).
Both age and time steps were chosen as $h_{\MA} = h_{\VE} = \tau = 1/12$.
Based on this selection, we linearly interpolated the data onto the grid.
Subsequently, Equations (\ref{EQUATION_FULL_DISCRETIZATION_ALGEBRAIC_1})--(\ref{EQUATION_FULL_DISCRETIZATION_ALGEBRAIC_3})
were solved using the Crank \& Nicholson method corresponding to selecting $\vartheta = 1/2$
and the output was back-transformed using the age-specific survival probabilities.
Finally, we restricted the simulation results onto the single-year-spaced grid.
Our \texttt{Matlab}-code can be downloaded from \texttt{MathWorks} under
\texttt{http://www.mathworks.com/matlabcentral/fileexchange/48072}

Table \ref{TABLE_US_POPULATION_2011} below
gives a comparison between the total male and female population in the U.S. as reported by \cite{USCB2013}
and as estimated from our simulation.
As Table \ref{TABLE_US_POPULATION_2011} suggests,
we underestimated both the male and female population by merely 2.54\% and 2.82\%, respectively.
Probably, this is due to the fact the immigration data are not sufficiently reliable
and tend to be somewhat underestimated in official surveys.
Though not being perfect, 
our estimate seem to outperform the expected precision of 4.1\% described in \cite{AbLo1997} for the decade 1970--1980.
Thus, our prediction seems to be rather accurate
even without accounting for the official marital status of population members unlike \cite{AbLo1997}.
\begin{table}[h!]
	\centering
	\begin{tabular}{ccccc}
		\toprule
		& \multicolumn{2}{c}{Total number} & \multicolumn{2}{c}{Relative error} \\
		& Men & Women & Men & Women \\
		\midrule
		Reported  & 153253317 & 158287949 & --     & -- \\
		Simulated & 149360262 & 153825899 & 2.54\% & 2.82\% \\
		\bottomrule
	\end{tabular}
	\caption{Summary on the U.S. population in 2011. \label{TABLE_US_POPULATION_2011}}
\end{table}

Table \ref{TABLE_US_POPULATION_2011_ERROR} gives the actual errors,
i.e., the discrepancy between the simulated and reported data in different norms.
Related to the total male or female population, the error never exceeded 3.68\% measured with respect to any $L^{p}$-norm, $p = 1, 2, \infty$.
\begin{table}[h!]
	\centering
	\begin{tabular}{ccccccc}
		\toprule
		& \multicolumn{2}{c}{$L^{1}$} & \multicolumn{2}{c}{$L^{2}$} & \multicolumn{2}{c}{$L^{\infty}$}\\
		         & Men     & Women   & Men    & Women  & Men    & Women \\
		\midrule
		Absolute & 5054906 & 5819685 & 702318 & 746205 & 218037 & 228500 \\
		Relative & 3.30\%  & 3.68\%  & 0.46\% & 0.47\% & 0.14\% & 0.14\% \\
		\bottomrule
	\end{tabular}
	\caption{Actual errors. \label{TABLE_US_POPULATION_2011_ERROR}}
\end{table}

\begin{figure}[h!]
	\centering
	\includegraphics[scale = 0.4]{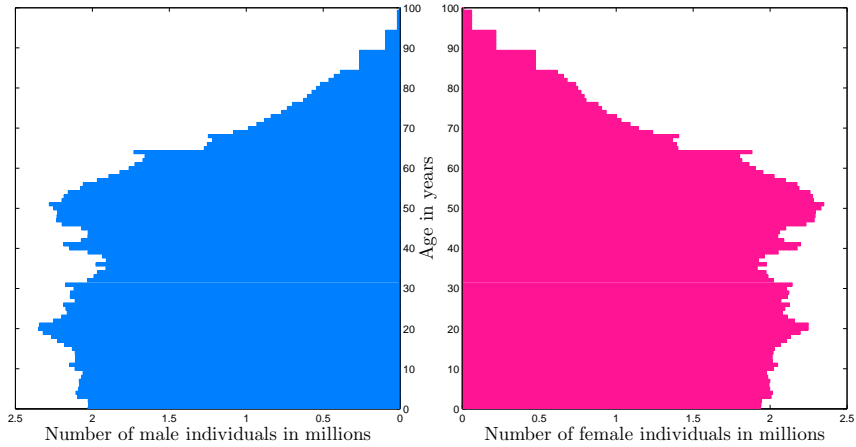}
	\caption{Reported age-sex-structure of the United States in 2011. \label{FIGURE_POPULATION_2011_REPORTED}}
\end{figure}

\begin{figure}[h!]
	\centering
	\includegraphics[scale = 0.4]{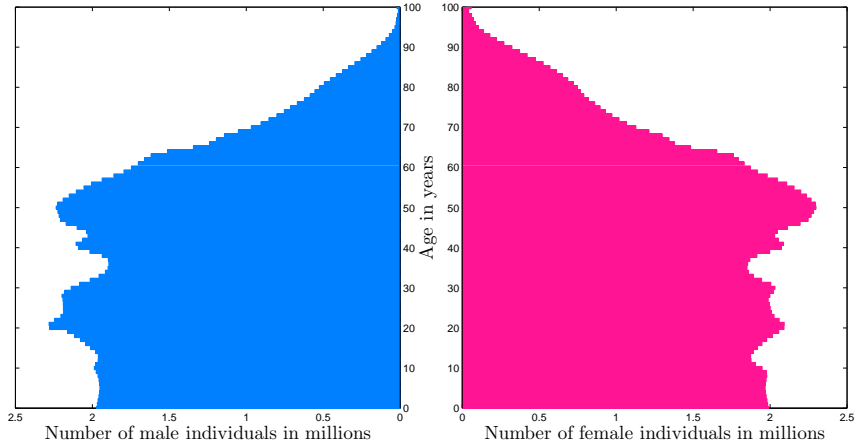}
	\caption{Simulated age-sex-structure of the United States in 2011. \label{FIGURE_POPULATION_2011_SIMULATED}}
\end{figure}

Finally, Figure \ref{FIGURE_POPULATION_2011_REPORTED} displays the U.S. population in 2011 as reported in \cite{USCB2013},
whereas Figure \ref{FIGURE_POPULATION_2011_SIMULATED} depicts the outcome of our numerical simulation for the same year.
Both Figures seem to be in a good accordance with each other
though the reported population looks somewhat ``spiky''.
Statistically, the latter can be explained by the fact
the data are binned and thus can exhibit such roughness patterns due to grouping (cf. \cite[Chapter 2]{Hae2004}).

\begin{appendix}
\section{Discrete Spaces and Operators}
Let $I = (a, b) \subset \RR$ be a bounded interval and let $X$ be a Hilbert space.
For $h > 0$ such that $N_{h} = \tfrac{b - a}{h} \in \mathbb{N}$, let $I$ be partitioned by an equidistant lattice $I^{h} = \{\xi^{h}_{k} \,|\, k = 0, \dots, N_{h}\}$
with $\xi^{h}_{k} = a + kh$, $k = 0, \dots, N_{h} = \tfrac{b - a}{h}$.
We define the discrete Lebesgue $L^{2}$-space
\begin{equation}
	L^{2}_{h}(I^{h}, X) := X^{I^{h}}, \quad
	\langle u^{h}, v^{h}\rangle_{L^{2}_{h}(I^{h}, X)} :=
	h \sum_{k = 0}^{N_{h}} \langle u_{h}(\xi^{h}_{k}), v^{h}(\xi^{h}_{k})\rangle_{X}
	\text{ for } u^{h}, v^{h} \colon I^{h} \to X. \notag
\end{equation}
For $X = \RR$, we simply write $L^{2}_{h}(I^{h})$.

Letting $\stackrel{\circ}{I^{h}} := \{\xi^{h}_{k} \,|\, k = 1, \dots, N_{h}\}$ and
$\stackrel{\circ}{\bar{I}^{h}} := \{\xi^{h}_{k} \,|\, k = 0, \dots, N_{h}-1\}$,
we define the backwards and forwards difference operators
\begin{equation}
	\begin{split}
		\partial^{\xi}_{h} \colon L^{2}_{h}(I^{h}, X) \to L^{2}_{h}(\stackrel{\circ}{I^{h}}, X), \quad
		\partial^{\xi}_{h} u := \xi^{h}_{k} \mapsto \tfrac{u(\xi^{h}_{k}) - u(\xi^{h}_{k-1})}{h} \text{ for } u \in L^{2}_{h}(I^{h}, X), \\
		\bar{\partial}^{\xi}_{h} \colon L^{2}_{h}(I^{h}, X) \to L^{2}_{h}(\stackrel{\circ}{\bar{I}^{h}}, X), \quad
		\bar{\partial}^{\xi}_{h} u := \xi^{h}_{k} \mapsto \tfrac{u(\xi^{h}_{k+1}) - u(\xi^{h}_{k})}{h} \text{ for } u \in L^{2}_{h}(I^{h}, X),
	\end{split}
	\notag
\end{equation}
respectively. Note that both $\partial^{h}_{\xi}$ and $\bar{\partial}^{h}_{\xi}$
are linear, bounded operators from $H^{1}(I, X)$ to $L^{2}_{h}(\stackrel{\circ}{\bar{I}^{h}}, X)$ and $L^{2}_{h}(\stackrel{\circ}{I^{h}}, X)$, respectively,
by the virtue of Sobolev embedding theorem.
We have the well-known summation by parts formula:
\begin{lemma}
	\label{LEMMA_SUMMATION_BY_PARTS}
	For $u \in L^{2}_{h}(I^{h}, \RR), v \in L^{2}_{h}(I^{h}, X)$, there holds
	\begin{equation}
		\sum_{\xi_{k} \in \stackrel{\circ}{I^{h}}} (\partial_{\xi}^{h} u)(\xi_{k}) v(\xi_{k}) =
		-\sum_{\xi_{k} \in \stackrel{\circ}{\bar{I}^{h}}} u(\xi_{k}) \bar{\partial}_{\xi}^{h} v(\xi_{k})
		+ u(b) v(b) - u(a) v(a). \notag
	\end{equation}
\end{lemma}

As an immediate consequence of \cite[Propositions 1.1.6 and 1.2.2]{ArBaHieNeu2001}, we have the following two lemmas
\begin{lemma}
	\label{LEMMA_APPENDIX_DIFFERENCE_OPERATOR_ESTIMATE}
	For any $u \in H^{1}(I, X)$, 
	there holds
	\begin{equation}
		\sum_{k = 1}^{N_{h}} \int_{\xi_{k-1}}^{\xi^{k}} \big\|\partial^{h}_{\xi} u^{h}(\xi_{k}) - \partial_{\xi} u(\xi)\big\|_{X}^{2} \mathrm{d}\xi \to 0 \text{ as } h \to 0. \notag
	\end{equation}
\end{lemma}

\begin{lemma}
	\label{LEMMA_APPENDIX_INTEGRATION_OPERATOR_ESTIMATE}
	Let $M \in C^{0}(\bar{I}, L(X))$. For any $u \in C^{0}(\bar{I}, X)$, there holds
	\begin{equation}
		\frac{1}{h} \Big\|\int_{I} M(\xi) u(\xi) \mathrm{d}\xi - h \sum_{k = 1}^{N_{h}} M(\xi^{h}_{k}) u(\xi^{h}_{k})\Big\|_{X} \to 0 \text{ as } h \to 0. \notag
	\end{equation}
\end{lemma}
\end{appendix}

\section*{Acknowledgments}
This work has been funded by a research grant from the Young Scholar Fund
supported by the Deutsche Forschungsgemeinschaft (ZUK 52/2) at the University of Konstanz, Konstanz, Germany.

\end{document}